\documentclass[10pt]{amsart}
\usepackage{amssymb}

\newtheorem{lemma}{Lemma}[section]

\newtheorem{theorem}[lemma]{Theorem}
\newtheorem{remark}[lemma]{Remark}
\newtheorem{proposition}[lemma]{Proposition}
\newtheorem{corollary}[lemma]{Corollary}
\newtheorem{definition}[lemma]{Definition}

\begin{document}

\title{Principal kinematic formulas for germs of closed definable sets
}

\author{Nicolas Dutertre}
\address{Laboratoire angevin de recherche en math\'ematiques, LAREMA, UMR6093, CNRS, UNIV. Angers, SFR MathStic, 2 Bd Lavoisier 49045 Angers Cedex 01, France.}
\email{nicolas.dutertre@univ-angers.fr}
 
\thanks{The author is partially supported by the ANR project LISA 17-CE400023-01}



\markboth{N. Dutertre}{Principal kinematic formulas for germs of closed definable sets}

\begin{abstract}
We prove two principal kinematic formulas for germs of closed definable sets in $\mathbb{R}^n$, that generalize the Cauchy-Crofton formula for the density due to Comte and the infinitesimal linear kinematic formula due to the author. In this setting, we do not integrate on the space of euclidian motions $SO(n) \ltimes \mathbb{R}^n$, but on the manifold $SO(n) \times \mathbb{S}^{n-1}$. 
\end{abstract}

\maketitle
\section{Introduction}
The search for kinematic formulas is one of the main goal of integral geometry. 
Such formulas have been proved in various contexts by various authors, for instance:
\begin{itemize}
\item[-] For convex bodies by Blaschke and Hadwiger (see \cite{KlainRota});
\item[-] For manifolds by Chern \cite{Chern} and manifolds with boundary by Santal\'o \cite{Santalo};
\item[-] For PL-sets by Cheeger, M\"uller and Schrader \cite{Cheeger};
\item[-] For sets with positive reach by Federer \cite{FedererTAMS,FedererBook} (see also \cite{RotherZahle});
\item[-] For subanalytic sets by Fu \cite{Fu94}, and more generally for sets definable in an o-minimal structure by Bernig, Br\"ocker and Kuppe \cite{BroeckerKuppe,BernigBroeckerMathNachr,BernigBroeckerFourier}.
\end{itemize}
There are many other situations where kinematic formulas hold,
but we cannot give here a complete list of all the interesting papers published on this topics.

In this paper, we are interested in the case of definable sets in an o-minimal structure. Definable sets are a generalization of semi-algebraic sets and globally subanalytic sets, we refer the reader to classical references \cite{DriesMiller,Dries,CosteOmin,Loi10,Rolin} for basic definitions and results on this topics. The study of the geometric properties of these objects was initiated by Fu \cite{Fu94}, who developed integral geometry for compact subanalytic sets. Using the technology of the normal cycle, he associated with every compact subanalytic set $X$ of $\mathbb{R}^n$ a sequence of curvature measures
$$\Lambda_0(X,-),\ldots,\Lambda_n(X,-),$$
called the Lipschitz-Killing measures, and he established several integral geometry formulas.
Among them, he proved the following  kinematic formulas:  for $k \in \{0,\ldots,n \}$, we have
$$\int_{SO(n) \ltimes \mathbb{R}^n} \Lambda_k(X \cap \gamma Y, U \cap \gamma V) d\gamma dx = \sum_{p+q=k+n} e(p,q,n)
\Lambda_p(X,U) \Lambda_q(X,V) ,$$
where $X$ and $Y$ are two compact subanalytic subsets of $\mathbb{R}^n$ and $U$ and $V$ are Borel subsets of $X$ and $Y$ respectively. We will state these formulas specifically in the next section.
In \cite{BroeckerKuppe} (see also \cite{BernigBroeckerMathNachr,BernigBroeckerFourier}), Br\"ocker and Kuppe gave a geometric characterization of these measures using stratified Morse theory, in the more general setting of definable sets.  

In \cite{Comte} Comte started the study of real equisingularity by proving that the density is continuous along the strata of a Verdier stratification of a subanalytic set (see also \cite{Valette}). The main tool to prove his result was a local Cauchy-Crofton formula for the density. He continued this work with Merle in \cite{ComteMerle} where a similar continuity result was established for the so-called local Lipschitz-Killing invariants (see also \cite{NguyenValette}). The tools for proving this continuity property are local linear kinematic formulas that generalize the Cauchy-Crofton formula for the density. 
These formulas will be explained in Section 3 but, roughly speaking, they relate the so-called polar invariants, which are mean-values of Euler characteristics of real Milnor fibres of generic projections, to the local Lipschitz-Killing invariants.

In \cite{DutertreJofSingProcTrot} we also established an infinitesimal linear kinematic formula. It is slightly different from the ones of Comte and Merle, because instead of using projections, we make ``infinitesimally small" translations of linear spaces. Let us recall it here because it is our main inspiration. We will use the following notations:
\begin{itemize}
\item $s_k$ is the volume of unit sphere $\mathbb{S}^k$ of dimension $k$ and $b_k$ is the volume of the unit ball $\mathbb{B}^k$ of dimension $k$,
\item for $k \in \{0,\ldots,n \}$, $G_{n}^k$ is the Grassmann manifold of $k$-dimension  linear spaces in $\mathbb{R}^{n}$ equipped with the $O(n)$-invariant Maurer-Cartan density (see for instance \cite{Santalo}, p.200),
$g_n^k$ is its volume,
\item if $P$ is a linear subspace of $\mathbb{R}^{n}$ of dimension $k$, $\mathbb{S}_P^{k-1}$ is the unit sphere in $P$,
\item in $\mathbb{R}^n$, $\mathbb{B}_\epsilon^n(x)$ is the closed ball of radius $\epsilon$ centered at $x$ and $\mathbb{S}_\epsilon^{n-1} (x)$ is the sphere of radius $\epsilon$ centered at $x$, if $x=0$, we simply write $\mathbb{B}^n_\epsilon$ and $\mathbb{S}^{n-1}_\epsilon$.
\end{itemize}
Let $(X,0) \subset (\mathbb{R}^n,0)$ be the germ of a closed definable set. We consider the following limits:
$$\Lambda_k^{\rm lim}(X,0):= \lim_{\epsilon \rightarrow 0} \frac{\Lambda_k(X,X\cap \mathbb{B}_\epsilon^n)}{b_k \epsilon^k}.$$
Let $H \in G_n^{n-k}$, $k \in \{1,\ldots, n \}$, and let $v$ be an element in $\mathbb{S}_{H^\perp}^{k-1}$. For $\delta >0$, we denote by $H_{v,\delta}$ the $(n-k)$-dimensional affine space $H+\delta v$ and we set 
$$\beta_0(H,v) = \lim_{\epsilon \rightarrow 0} \lim_{\delta \rightarrow 0} \Lambda_0(H_{\delta,v} \cap X, H_{\delta,v} \cap X \cap \mathbb{B}_\epsilon^n),$$
and  
$$\beta_0 (H) =\frac{1}{s_{k-1}}\int_{\mathbb{S}_{H^\perp}^{k-1}} \beta_0 (H,v) dv.$$ In \cite{DutertreJofSingProcTrot} Theorem 5.5, we proved that for $k \in \{1,\ldots,n \}$ 
$$ \Lambda_k^{\rm lim}(X,0) = \frac{1}{g_n^{n-k}} \int_{G_n^{n-k}} \beta_0 (H) dH.$$

In view of this formula and since it is possible to make ``infinitesimally small" translations of any definable set, the question that motivated us was the following: Is it possible to establish a kinematic formula for germs of closed definable sets or, in other words, can we replace the $(n-k)$-plane $H$ with any germ of closed definable set? The goal of this paper is to provide a positive answer to this question.

Let us present the main results of the paper. Let $(X,0) \subset (\mathbb{R}^n,0)$ be the germ of a closed definable set. To such a germ, we associate two sequences of real numbers: the polar invariants $\sigma_k(X,0)$, $k=0,\ldots,n$, and the above limits
$\Lambda_k^{\rm lim}(X,0)$.
Let $(Y,0) \subset (\mathbb{R}^n,0)$ be another germ of closed definable set 
and let 
$$\sigma(X,Y,0) = \frac{1}{s_{n-1}^2} 
\int_{SO(n) \times \mathbb{S}^{n-1}}  \lim_{\epsilon \to 0} \lim_{\delta \to 0^+}
 \chi \left( X \cap (\gamma Y +\delta v) \cap \mathbb{B}_\epsilon^n \right) d\gamma dv.$$
Here $SO(n)$ is equipped with the Haar measure $d\gamma$, normalized in such a way that the volume of $SO(n)$ is $s_{n-1}$, $\mathbb{S}^{n-1}$ is equipped with the usual Riemanniann measure (or density) $dv$ and $SO(n) \times \mathbb{S}^{n-1}$ with the product measure $d\gamma dv$. 
Our first infinitesimal principal kinematic formula takes the following form (Theorem \ref{Theorem8.15}):
$$ \sigma(X,Y,0)= \sum_{i=0}^n \Lambda_i^{\rm lim}(X,0) \cdot \sigma_{n-i}(Y,0).$$
When $X$ and $Y$ have complementary dimensions, this gives a Bezout type formula, since the integrand of the left-hand side is a number of intersection points and the right-hand side is the product of the densities of $X$ and $Y$ at the origin (see Corollary \ref{CorollaryofTh8.15}).
Then we set  
$$\displaylines{
\quad \Lambda_0^{{\rm lim}}(X,Y,0) \hfill \cr
\hfill = \frac{1}{s_{n-1}^2} 
\int_{SO(n) \times \mathbb{S}^{n-1}}  \lim_{\epsilon \to 0} \lim_{\delta \to 0^+}
\Lambda_0 \left(X \cap (\gamma Y + \delta v) ,  X \cap (\gamma Y + \delta v) \cap \mathbb{B}_\epsilon^n \right) d\gamma dv. \cr
}$$

In Theorem \ref{Theorem8.16}, we establish our second infinitesimal principal kinematic formula:
 $$\Lambda_0^{\rm lim} (X,Y,0)= \sum_{i=0}^n \Lambda_i^{\rm lim}(X,0) \cdot \Lambda_{n-i}^{\rm lim}(Y,0).$$
This formula is a corollary of Theorem \ref{Theorem8.15} and the Gauss-Bonnet formula proved in \cite{DutertreJofSingProcTrot}. 

Throughout the paper, we will also use the following notations and conventions:
\begin{itemize}
\item for $v \in \mathbb{R}^n$, the function $v^* : \mathbb{R}^n \rightarrow \mathbb{R}$ is defined by $v^*(y)= \langle v, y \rangle$,
\item for $x \in \mathbb{R}^n$, $\vert x \vert$ denotes the usual Euclidean norm,
\item if $X \subset \mathbb{R}^{n}$,   $\overline{X}$ is its
topological closure, $\mathring{X}$ its topological interior,
\item  when it makes sense, ${\rm vol}(X)$ means the volume of the set $X$ and $\chi_c(X)$ its Euler characteristic for Borel-Moore homology.
\end{itemize}

The paper is organized as follows. In Section 2, we recall the notion of  stratified critical points and the definition of the Lipschitz-Killing measures. We also state kinematic formulas. In Section 3, we recall the Gauss-Bonnet formula for real Milnor fibres proved by the author in \cite{DutertreJofSingProcTrot}, and the infinitesimal linear kinematic formulas proved by Comte \cite{Comte}, Comte and Merle \cite{ComteMerle} and the author \cite{DutertreJofSingProcTrot}. Section 4 contains several topological and geometrical lemmas that will be useful in the next sections. In Section 5, we prove a new spherical kinematic formula for definable sets. Combining this formula with Hardt's theorem \cite{Hardt,CosteOmin,DriesMiller}, we obtain a new kinematic formula for definable subsets of the unit ball in Section 6. We apply this formula in Section 7 to get our first principal kinematic formula for closed conic definable sets. In Section 8, we prove our first principal kinematic formula in the general case using the previous case and tangent cones, and then our second principal kinematic formula. Finally Section 9 contains two other kinematic formulas.

\section{Stratified critical points and Lipschitz-Killing curvatures}

\subsection{Stratified critical points}
Let $X \subset \mathbb{R}^n$ be a compact definable set equipped with a finite definable Whitney stratification $\mathcal{S}=\{ S_a \}_{a \in A}$. The fact that such a stratification exists is due to Loi \cite{Loi98} (see also \cite{NguyenTrivediTrotman}).

Let $f: X \to \mathbb{R}$ be a definable function. We assume that $f$ is the restriction to $X$ of a $C^2$ definable function $F : U \to \mathbb{R}$, where $U$ is an open neighborhood of $X$ in $\mathbb{R}^n$. A point $p$ in $X$ is a (stratified) critical point of $f$ if $p$ is a critical point of $f_{\vert S}$, where $S$ is the stratum that contains $p$.

\begin{definition}\label{Definition2.1}
{\rm Let $p \in X$ be an isolated critical point of  $f : X \rightarrow \mathbb{R}$. 
The index of $f$ at $p$ is defined by
$$\hbox{ind}(f,X,p)=1-\chi \left( X \cap \{f=f(p_i)-\delta \} \cap \mathbb{B}_\varepsilon^n    (p) \right),$$
where $0< \delta \ll \varepsilon \ll 1$. If $p \in X$ is not a critical point of $f$, we set ${\rm ind}(f,X,p)=0$.} 
\end{definition}
Since we are in the definable setting, this index is well-defined thanks to Hardt's  theorem \cite{Hardt,CosteOmin,DriesMiller}. 

\begin{theorem}\label{PoincareHopf}
Assume that $f : X \to \mathbb{R}$ has a finite number of critical points $\{p_1,\ldots,p_s\}$. Then the following equality holds:
$$\chi (X)= \sum_{i=1}^s {\rm ind} (f,X,p_i).$$
\end{theorem}
\proof See Theorem 3.1 in \cite{DutertreManuscripta12}. When $f$ is a Morse stratified function, this follows from \cite{GoreskyMacPherson}. \endproof

\subsection{Lipchitz-Killing curvatures}
In this subsection, we present the Lipschitz-Killing measures of a definable set in an o-minimal structure.  We describe Br\"ocker and Kuppe's approach \cite{BroeckerKuppe}.

Let $X \subset \mathbb{R}^n$ be a compact definable set equipped with a finite definable Whitney stratification $\mathcal{S}=\{ S_a \}_{a \in A}$.  

Let  us fix a stratum $S$.
For $k \in \{0,\ldots,d_S\}$, $d_S= {\rm dim} S$, let $\lambda_k^S : S \rightarrow \mathbb{R}$ be defined by
$$ \lambda_k^S(x) = \frac{1}{s_{n-k-1}} \int_{\mathbb{S}_{T_x S^\perp} }{\rm ind}_{\rm nor}(v^*,X,x) \sigma_{d_S-k} (II_{x,v}) dv,$$
where $II_{x,v}$ is the second fundamental form on $S$ in the direction of $v$ and where $\sigma_{d_S-k} (II_{x,v})$ is the $(d_S-k)$-th elementary symmetric function of its eigenvalues. The index ${\rm ind}_{\rm nor}(v^*,X,x)$ is defined as follows:
$${\rm ind}_{\rm nor}(v^*,X,x)= 1-\chi \Big( X \cap N_x \cap \mathbb{B}_\epsilon^n(x) \cap \{ v^*= v^*(x)-\delta \} \Big),$$
where $0 < \delta \ll \epsilon \ll 1$ and $N_x$ is a normal (definable) slice to $S$ at $x$ in $\mathbb{R}^n$ . 
When $v^*_{\vert X}$ has a stratified Morse critical point at $x$, it coincides with the normal Morse index at $x$ of a function $f : \mathbb{R}^n \rightarrow \mathbb{R}$ such that $f_{\vert X}$ has a stratified Morse critical point at $x$ and $\nabla f (x) = v$. For $k \in \{d_S+1,\ldots,n\}$, we set $\lambda_k^S(x)=0$. 

If $S$ has dimension $n$ then for all $x \in S$, we put $\lambda_0^S(x)=\cdots=\lambda_{n-1}^S(x)=0$ and $\lambda_n^S(x)=1$. If $S$ has dimension $0$ then ${\rm ind}_{\rm nor}(v^*,X,x)= {\rm ind}(v^*,X,x)$ and we set
$$\lambda_0^S(x)= \frac{1}{s_{n-1}} \int_{\mathbb{S}^{n-1}}  {\rm ind}(v^*,X,x) dv,$$
and $\lambda_k^S(x)=0$ if $k>0$. 
\begin{definition}
{\rm For every Borel set $U \subset X$ and for every $k \in \{0,\ldots,n\}$, we define $\Lambda_k(X,U)$ by
$$\Lambda_k(X,U)= \sum_{a \in A} \int_{S_a \cap U} \lambda_k^{S_a} (x) dx.$$}
\end{definition}
These measures $\Lambda_k(X,-)$ are called the Lipschitz-Killing measures of $X$. Note that for any Borel set $U$ of $X$, we have
$$\Lambda_{d+1}(X,U)= \cdots=\Lambda_n(X,U)=0,$$ 
and $\Lambda_{d}(X,U)= \mathcal{H}_{d}(U)$, where $d$ is the dimension of $X$ and $\mathcal{H}_{d}$ is the $d$-th dimensional Hausdorff measure in $\mathbb{R}^n$. 
If $X$ is smooth then for $k \in \{0,\ldots,d \}$, $\Lambda_k(X,U)$ is equal to 
$$\frac{1}{s_{n-k-1} } \int_U K_{d-k}(x) dx.,$$
where $K_{d-k}$ denotes the $(d-k)$-th Lipschitz-Killing curvature.

As in the smooth case, the measure $\Lambda_0(X,-)$ satisfies an exchange formula (see \cite{BroeckerKuppe}).
\begin{proposition}\label{ExchangeFormula}
For every Borel set $U \subset X$, we have 
$$\Lambda_0(X,U) =\frac{1}{s_{n-1}} \int_{S^{n-1}} \sum_{x \in U}  {\rm ind} (v^*,X,x) dv.$$
\end{proposition}
For $U=X$ and by Theorem \ref{PoincareHopf}, we see that a special case of this exchange formula is the Gauss-Bonnet formula $\Lambda_0(X,X)=\chi(X)$.

The Lipschitz-Killing measures satisfy the kinematic formula (see \cite{Fu94,BroeckerKuppe,BernigBroeckerMathNachr,BernigBroeckerFourier}).                                                              
We provide the group $SO(n) \ltimes \mathbb{R}^n$ of all euclideans motions of $\mathbb{R}^n$ with the product measure $d\gamma dx$, where the canonical Haar measure $d\gamma$ is normalized such that ${\rm vol}\left(SO(n)\right)=1$.
\begin{proposition}\label{KinematicFormula}
Let $X\subset \mathbb{R}^n$ and $Y\subset \mathbb{R}^n$ be two compact definable sets and let $U \subset X$ and $V \subset Y$ be two Borel sets.  For $k \in \{0,\ldots,n \}$, the following kinematic formula holds:
$$\int_{SO(n) \ltimes \mathbb{R}^n} \Lambda_k(X \cap \gamma Y, U \cap \gamma V) d\gamma dx = \sum_{p+q=k+n} e(p,q,n)
\Lambda_p(X,U) \Lambda_q(X,V) ,$$
where $e(p,q,n)=\frac{s_{p+q-n} s_n}{s_p s_q}$.
\end{proposition}
For $k=0$, the above formula is called the principal kinematic formula.
A particular case of the kinematic formula is the linear kinematic formula. Let $A_{n}^k$ be the affine grassmannian of $k$-dimensional affine spaces in $\mathbb{R}^{n}$. It is a fibre bundle over $G_n ^k$ with fibre $\mathbb{R}^{n-k}$. We equip $A_n^k$ with the product measure denoted by $dE$.
\begin{proposition}\label{LinearKinematicFormula}
Let $X\subset \mathbb{R}^n$ be a compact definabet set and let $U \subset X$ be a Borel set. For $k \in \{0,\ldots,n \}$, we have
$$\Lambda_{n-k}(X,U) = \frac{1}{g_n^k} \cdot \frac{1}{e(k,n-k,n)}\int_{A_n^k} \Lambda_0(X \cap E, X \cap E \cap U) dE.$$
\end{proposition}

In Section 5, we will consider definable subsets of the unit sphere $\mathbb{S}^{n-1}$. For such sets, one can define spherical Lipschitz-Killing measures. These measure are defined in \cite{BernigBroeckerFourier}, Section 3 (see also \cite{DutertreGeoDedicata2012}). Their definition is very similar to the definition of the above Lipschitz-Killing measures. For $X \subset \mathbb{S}^{n-1}$ and $k \in \{0,\ldots,n-1\}$, we will denote by $\tilde{\Lambda}_k(X,-)$ the $k$-th spherical Lispchitz-Killing measures. The spherical Lipschitz-Killing measures satisfy a Gauss-Bonnet formula (\cite{BernigBroeckerFourier}, Theorem 1.2) and a spherical kinematic formula (\cite{BernigBroeckerFourier}, Theorem 4.4).

\section{Some topological ang geometrical properties of definable sets}
In this section, we review some results on the local topology and geometry of closed definable sets. Let $(X,0)$ be the germ of a closed definable set. For convenience, we will work with a small representative that we denote by $X$ as well. We assume that this representative is included in a an open bounded neighborhood $U$ of $0$. 

\subsection{The Gauss-Bonnet formula for real Milnor fibres}

We can equip $X$ with a finite Whitney stratification $\mathcal{S}=\{S_\alpha\}_{\alpha \in A}$ such that $0 \in \overline{S_\alpha}$ (this is possible taking a smaller representative if necessary). 

Let $\rho_i : U \to \mathbb{R}$, $i=1,2$, be two continuous definable functions of class $C^2$ on $U\setminus \{0\}$, such that $\rho_i^{-1}(0)=\{0\}$ and $\rho_i (x) \ge 0$ for all $x\in X$. It is well-known that there exists $\epsilon_i >0$ such that for $0< \epsilon \le \epsilon_i$, $\rho_i^{-1}(\epsilon)$ intersects $X$ transversally in the stratified sense (see \cite{DutertreJofSingProcTrot} Lemma 2.1), and that the topological type of $\rho_i^{-1}(\epsilon) \cap X$ does not depend on $\epsilon$. Moreover, as explained by Durfee in \cite{Durfee}, Lemma 1.8 and Corollary 3.6, there is a neighborhood $\Omega$ of $0$ in $\mathbb{R}^n$ such that for every stratum $S$ of $X$, $\nabla ({\rho_1}_{\vert S})$ and $\nabla ({\rho_2}_{\vert S})$ do not point in opposite direction in $\Omega \setminus \{0\}$. Applying Durfee's argument (\cite{Durfee}, Proposition 1.7 and Proposition 3.5), we see that $\rho_1^{-1}(\epsilon) \cap X$, $0< \epsilon \le \epsilon_1$, and $\rho_2^{-1}(\epsilon') \cap X$, $0< \epsilon' \le \epsilon_2$, are homeomorphic. The link of $X$ at $0$, denoted by ${\rm Lk}(X)$, is the set $X \cap \rho^{-1}(\epsilon)$, $0 < \epsilon \ll 1$, where $\rho: U \to \mathbb{R}$ is a continuous definable function of class $C^2$ on $U\setminus \{0\}$, such that $\rho^{-1}(0)=\{0\}$ and $\rho(x) \ge 0$ for all $x\in X$. We will call such a function $\rho$ a distance function to the origin. By the above discussion, the topological type of ${\rm Lk}(X)$ does not depend on the choice of the definable distance function to the origin (actually to define the link, we do not need to assume that $\rho$ is $C^2$ on $U \setminus \{0\}$, continuity is enough). 

Let $f: (X,0) \to (\mathbb{R},0)$ be the germ of a definable function. We assume that $f$ is the restriction to $X$ of a $C^2$ definable function $F : U \to \mathbb{R}$. We denote by $X^f$ the set $f^{-1}(0)$ and by \cite{Bekka,TaLeLoi}, we can equip $X$ with a definable Thom stratification $\mathcal{V}=\{V_\beta\}_{\beta \in B}$ adapted to $X^f$. This means that $\{ V_\beta \ \vert \ V_\beta \nsubseteq X^f\}$ is a Whitney stratification of $X \setminus X^f$ and that for any pair of strata $(V_\beta,V_{\beta'})$ with $V_\beta \nsubseteq X^f$ and $V_{\beta'} \subset X^f$, the Thom condition is satisfied. 

Note that if $f : (X,0) \to (\mathbb{R},0)$ has an isolated stratified critical point at $0$, where $X$ is equipped with the above Whitney stratification $\mathcal{S}=\{S_\alpha\}_{\alpha \in A}$, then the following stratification:
$$\left\{ S_\alpha \setminus f^{-1}(0), S_\alpha \cap (f^{-1}(0) \setminus \{0\}), \{0\} \ \vert \ \alpha \in A \right\},$$
is a Thom stratification of $X$ adapted to $X^f$. 

As explained above, there is $\epsilon'_1 >0$ such that for $0< \epsilon \le \epsilon'_1$, $\rho_1^{-1} (\epsilon) $ intersects $X^f$ transversally. The Thom condition implies that there exists $\delta_\epsilon >0$ such that for each $\delta$ with $0< \delta \le \delta_\epsilon$, $\rho_1^{-1}(\epsilon)$ intersects $f^{-1}(\delta)$ transversally as well. Hence the set $f^{-1}(\delta) \cap \{ \rho_1 \le \epsilon \}$ is a Whitney stratified set equipped with the following stratification:
$$\left\{ f^{-1}(\delta) \cap V_\beta \cap \{ \rho_1 < \epsilon \},
f^{-1}(\delta) \cap V_\beta \cap \{ \rho_1 = \epsilon \} \ \vert \ V_\beta \nsubseteq X^f \right\}.$$
Moreover, taking $\epsilon'_1$ and $\delta_\epsilon$ smaller if necessary, the topological types of $f^{-1}(\delta) \cap \{ \rho_1 \le \epsilon \}$ and $f^{-1}(\delta)\cap \{ \rho_1 = \epsilon \}$ do not depend on the couple $(\epsilon,\delta)$. To see this, it is enough to adapt the proof of Lemma 2.1 in \cite{DutertreJofSingNonIsol} to the stratified case. The same fact is true for negative values of $f$.

Of course, we can make the same construction with $\rho_2$ instead of $\rho_1$. But as above, there is a neighborhood $\Omega'$ of $0$ in $\mathbb{R}^n$ such that for every stratum $W$ of $X^f$, $\nabla ({\rho_1}_{\vert W})$ and $\nabla ({\rho_2}_{\vert W})$ do not point in opposite direction. Let us choose $\epsilon'>0$ and $\epsilon >0$ such that $\{ \rho_2 \le \epsilon'\} \subsetneq \{ \rho_1 \le \epsilon \} \subset \Omega'$. If $\epsilon$, $\epsilon'$ and $\delta$ are sufficiently small then, for every stratum $V \nsubseteq X^f$, $\nabla ({\rho_1}_{\vert V \cap f^{-1}(\delta)})$ and $\nabla ({\rho_2}_{\vert V \cap f^{-1}(\delta)})$ do not point in opposite direction in $\{ \rho_1 \le \epsilon\} \setminus \{ \rho_2 < \epsilon'\}$. Otherwise, by Thom $(a_f)$-condition, we would find a point $p$ in $X^f \cap (\{ \rho_1 \le \epsilon\} \setminus \{ \rho_2 < \epsilon'\})$ such that either $\nabla ({\rho_1}_{\vert W}) (p)$ or $\nabla ({\rho_2}_{\vert W}) (p)$ vanish or $\nabla ({\rho_1}_{\vert W}) (p)$ and $\nabla ({\rho_2}_{\vert W}) (p)$ point in opposite direction, where $W$ is the stratum of $X^f$ that contains $p$ (see the proof of Lemma 3.7 in \cite{DutertreJofSingProcTrot}). This is impossible if we are sufficiently close to the origin. Applying Durfee's argument mentioned above, we see that $f^{-1}(\delta) \cap \{\rho_1 \le \epsilon \}$ is homeomorphic to $f^{-1}(\delta) \cap \{\rho_2 \le \epsilon' \}$ and that $f^{-1}(\delta) \cap \{\rho_1 = \epsilon \}$ is homeomorphic to $f^{-1}(\delta) \cap \{\rho_2 = \epsilon' \}$. 

The positive (resp. negative) Milnor fibre of $f$ is the set $f^{-1}(\delta) \cap \{ \rho \le \epsilon \}$ (resp. $f^{-1}(-\delta) \cap \{ \rho \le \epsilon \}$), where $0< \delta \ll \epsilon \ll 1$ and $\rho$ is a distance function to the origin. The set $f^{-1}(\pm \delta) \cap \{ \rho = \epsilon \}$ is the boundary of the Milnor fibre.
By the previous discussion, the topological type of the positive (resp. negative) Milnor fibre and the topological type of its boundary do not depend on the choice of the definable distance function to the origin. 

In \cite{DutertreJofSingProcTrot}, we considered a second definable function-germ $g : (\mathbb{R}^n,0) \to (\mathbb{R},0)$ and we assumed that $g$ was the restriction to $X$ of a $C^2$ definable function $G : U \to \mathbb{R}$. Moreover, we assumed that $g$ satisfied the following two conditions:
\begin{itemize}
\item Condition (A): $g : (X,0) \to (\mathbb{R},0)$ has an isolated critical point at $0$.
\item Condition (B): the relative polar set 
$$\Gamma_{f,g}= \sqcup_{V_\beta \nsubseteq X^f} \Gamma_{f,g}^{V_\beta}=
\sqcup_{V_\beta \nsubseteq X^f} \left\{ x \in V_\beta \ \vert \ {\rm rank} \left[ \nabla (f_{\vert V_\beta}) (x), \nabla (g_{\vert V_\beta})(x) \right] < 2 \right\}$$
is a $1$-dimensional $C^1$ definable set (possibly empty) in a neighborhood of the origin.
\end{itemize}
We wrote $\Gamma_{f,g}= \sqcup_{i=1}^l \mathcal{B}_i$, where each $\mathcal{B}_i$ is a definable connected curve, and we considered the intersections points of $\Gamma_{f,g}$ with $f^{-1}(\delta) \cap \mathbb{B}_\epsilon^n$: 
$$\Gamma_{f,g} \cap (f^{-1}(\delta) \cap \mathbb{B}_\epsilon^n ) = \sqcup_{i=1}^l \mathcal{B}_i \cap (f^{-1}(\delta) \cap \mathbb{B}_\epsilon^n) =\left\{ p_1^{\delta,\epsilon},\ldots,p_r^{\delta,\epsilon} \right\},$$
where $0 < \vert \delta \vert \ll \epsilon \ll 1$. The points $p_i^{\delta,\epsilon}$ are exactly the critical points of $g$ on $f^{-1}(\delta) \cap \mathring{\mathbb{B}_\epsilon^n}$. Then we set 
$$ I(\delta,\epsilon, g) =\sum_{i=1}^r {\rm ind}(g,f^{-1}(\delta),p_i^{\delta,\epsilon}),$$
$$ I(\delta,\epsilon, -g) =\sum_{i=1}^r {\rm ind}(-g,f^{-1}(\delta),p_i^{\delta,\epsilon}),$$
and in \cite{DutertreJofSingProcTrot}, Theorem 3.10, we related $I(\delta,\epsilon, g) +I(\delta,\epsilon, -g)$, with $0 < \vert \delta \vert \ll \epsilon \ll 1$, to the topology of the Milnor fibre and its boundary.

Let us give now a new characterization of $I(\delta,\epsilon,g)$ and $I(\delta,\epsilon,-g)$ independent on $\delta$ and $\epsilon$. Let us fix a connected component $\mathcal{B}$ of $\Gamma_{f,g}$. We can assume that $f$ is stricly increasing on $\mathcal{B}$ and we put $\mathcal{B} \cap f^{-1}(\delta)= \{p^\delta \}$ for $\delta >0$.
\begin{lemma}\label{Lemma3.1}
There exists $\delta_0 >0$ such that for $0< \delta \le \delta_0$, the function 
$\delta \mapsto {\rm ind}(g,f^{-1}(\delta),p^\delta)$ is constant on $]0,\delta_0]$.
\end{lemma}
\proof Let $d : \mathbb{R}^n \to \mathbb{R}$ be the distance function to $\mathcal{B}$. It is a continuous definable function on an open definable neighborhood $\mathcal{O}$ of $\mathcal{B}$. Let 
$$A = \left\{ x \in X \ \vert \ \exists p \in \mathcal{B} \hbox{ such that } 
f(x)=f(p) \hbox{ and } g(x) \le g(p) \right\}.$$
It is a definable subset of $X$. Let $\rho : A \cap (\mathcal{O} \setminus \mathcal{B}) \to \mathbb{R}^2$ be the mapping defined by $\rho(x)=(f(x),d(x))$. By Hardt's theorem \cite{Hardt,CosteOmin,DriesMiller}, there is a partition of $]0,+\infty[ \times ]0,+\infty[$ into finitely many definable sets such that $\rho$ is trivial over each of this set. Let us denote by $\Delta$ the union of the members of this partition which have dimension less than or equal to 1. By Hardt's theorem again, the set $$\left\{ \nu \in ]0,+\infty[ \ \vert \ \Delta \cap (\{\nu\} \times ]0,+\infty[ ) \hbox{ has dimension } 1 \right\}$$ is finite. For $\nu >0$, the function $r(\nu) = {\rm inf} \{ \epsilon' \ \vert \ (\nu,\epsilon') \in \Delta \}$ is definable and by the previous remark, there is $\nu_1 >0$ such that $r(\nu) >0$ for $0< \nu < \nu_1$. Hence by the Monotonicity Theorem (see \cite{CosteOmin}, Theorem 2.1 or \cite{DriesMiller}, 4.1), there is $0< \delta_0 < \nu_1$ such that $r$ is continuous, monotone and strictly positive on $]0,\delta_0]$. Moreover the function
$(\delta,\epsilon') \mapsto \chi \left( A \cap \{f=\delta\} \cap \{d=\epsilon'\} \right)$ is constant on $\{ (\delta,\epsilon') \ \vert   \ 0< \delta < \delta_0, 0< \epsilon' < r(\delta) \}$.
But, by Lemma 3.1 in \cite{DutertreJofSingNonIsol} and the above discussion on the topology of the link, we have
$$\displaylines{
\qquad \qquad {\rm ind}(g,f^{-1}(\delta),p^{\delta}) = 1 - \chi \left( \{g \le g(p^{\delta})\} \cap \{f=\delta\} \cap \{d=\epsilon'\} \right) \hfill \cr
\hfill = 1-\chi \left( A \cap \{f=\delta\} \cap \{d=\epsilon'\} \right). \qquad \qquad
}$$ 
We conclude that the function $\delta \mapsto {\rm ind}(g,f^{-1}(\delta),p^\delta)$ is constant on $]0,\delta_0]$.\endproof
Of course, a similar result holds for negative values of $f$. 

By the general \L ojasiewicz inequality (see \cite{BroeckerKuppe}, Corollary 1.5.2), there exists a continuous definable function $\psi : (\mathbb{R},0) \to (\mathbb{R},0)$ such that $\vert p \vert \le \psi (f(p))$ for $p \in \overline{\mathcal{B}}$. Moreover $\psi$ is of class $C^2$ in an open neighborhood of $0$ and $\psi(u)>0$ for $u>0$. Let us fix $\epsilon >0$ small. If $0< \delta < \psi^{-1} (\frac{\epsilon}{4})$ then $\vert p \vert \le \frac{\epsilon}{4}$ for $p \in \mathcal{B} \cap f^{-1}(\delta)$. 

Since $\Gamma_{f,g}$ consists of a finite number of branches, we can conclude that for $\epsilon >0$ there exists $\delta_0>0$ such that for $0 < \vert \delta \vert \le \delta_0$, $\Gamma_{f,g} \cap f^{-1}(\delta) \subset \mathbb{B}_{\frac{\epsilon}{4}}^n$, and so 
$\Gamma_{f,g} \cap (f^{-1}(\delta) \cap \mathbb{B}_\epsilon^n ) = \Gamma_{f,g} \cap f^{-1}(\delta)$. With the above notation, this means that $p_i^{\delta,\epsilon}=p_i^\delta$ for $0< \vert \delta \vert \ll \epsilon \ll 1$ and $i \in \{1,\ldots,r\}$. For $i \in \{1,\ldots,r\}$, let $\tau_i(g)$ (resp. $\tau_i(-g)$) be the value that the function $\delta \mapsto {\rm ind}(g,f^{-1}(\delta),p_i^\delta)$ (resp. ${\rm ind}(-g,f^{-1}(\delta),p_i^\delta)$) takes close to the origin. We deduce the following relations:
$$\lim_{\epsilon \to 0} \lim_{\delta \to 0^+} 
I(\delta,\epsilon, g) +I(\delta,\epsilon, -g) = \sum_{i \ \vert \ f>0 \ \hbox{\tiny on } \mathcal{B}_i} \tau_i(g)+ \tau_i (-g),$$
$$\lim_{\epsilon \to 0} \lim_{\delta \to 0^-} 
I(\delta,\epsilon, g) +I(\delta,\epsilon, -g) = \sum_{i \ \vert \ f<0 \ \hbox{\tiny on } \mathcal{B}_i} \tau_i(g)+ \tau_i (-g).$$
Of course, the same study can be done with another definable distance function to the origin and so, the two limits $\lim_{\epsilon \to 0} \lim_{\delta \to 0^+} 
I(\delta,\epsilon, g) +I(\delta,\epsilon, -g)$ and $\lim_{\epsilon \to 0} \lim_{\delta \to 0^-} 
I(\delta,\epsilon, g) +I(\delta,\epsilon, -g)$ do not depend on the distance function to the origin chosen to define the Milnor fibre of $f$. 

Applying this study to a linear form $v^*$, where $v$ is generic in $\mathbb{S}^{n-1}$, we established in \cite{DutertreJofSingProcTrot}, Theorem 4.5, an infinitesimal Gauss-Bonnet formula for the Milnor fibre of $f$. We will use only this formula for functions with an isolated stratified critical point at $0$. 
Namely if $X$ is equipped with a Whitney stratification for which $f:(X,0) \to (\mathbb{R},0)$ has an isolated stratified critical point at $0$, then (\cite{DutertreJofSingProcTrot}, Corollary 4.6)
$$\displaylines{
\qquad \lim_{\epsilon \to 0} \lim_{\delta \to \pm 0} 
\Lambda_0 \left(f^{-1}(\delta), f^{-1}(\delta) \cap \mathbb{B}_\epsilon^n \right)= 
\lim_{\epsilon \to 0} \lim_{\delta \to \pm 0}  \chi \left( f^{-1}(\delta) \cap \mathbb{B}_\epsilon^n \right) \hfill \cr
\hfill -\frac{1}{2} \chi \left( {\rm Lk} (X^f) \right)
-\frac{1}{2 s_{n-1}}\int_{\mathbb{S}^{n-1}} \chi \left({\rm Lk} (X^f  \cap \{v^*=0\}) \right)dv. \qquad \cr
}$$
The proof of this Gauss-Bonnet formula relies on the following exchange formula:
$$\displaylines{
\qquad \lim_{\epsilon \to 0} \lim_{\delta \to \pm 0} 
\Lambda_0 \left(f^{-1}(\delta), f^{-1}(\delta) \cap \mathbb{B}_\epsilon^n \right) \hfill \cr
\hfill =  \frac{1}{2 s_{n-1}}\int_{\mathbb{S}^{n-1}} \lim_{\epsilon \to 0} \lim_{\delta \to \pm 0}  \big[ 
I(\delta,\epsilon,v^*) + I(\delta,\epsilon,-v^*) \big] dv. \qquad \cr
}$$
But we have explained above that $\lim_{\epsilon \to 0} \lim_{\delta \to 0}  
I(\delta,\epsilon,v^*) + I(\delta,\epsilon,-v^*) $ does not depend on the choice of the distance function to the origin used to define the Milnor fibre of $f$. Therefore the relations proved in \cite{DutertreJofSingProcTrot}, Theorem 4.5 and Corollary 4.6, are also valid if we replace the usual euclidian distance function by any definable distance function to the origin. This remark will be important in the  next sections.

\subsection{Linear kinematic formulas for germs of closed definable sets} Let us recall the definition of the polar invariants \cite{ComteMerle}. 
Let $k \in \{1,\ldots,n\}$ and let $P \in G_n^k$. We denote by $\pi_P : X \to P$ the orthogonal projection on $P$. For $P$ generic in $G_n^k$, Comte and Merle established the existence of an open and dense definable germ $(K^P,0) \subset (P,0)$ such that, if $K^P= \cup_{i=1}^{N_P} K_i^P$ denotes its decomposition into connected components, then the function 
$$K_i^P \mapsto \chi_i^P := \lim_{\epsilon \to 0} \lim_{y \in K_i^P \atop y \to 0}
\chi \left( \pi_P^{-1}(y) \cap X \cap \mathbb{B}_\epsilon^n \right)$$
is well-defined. Then they set the following definition:
\begin{definition}\label{Definition3.2}
{\rm Let $k \in \{1,\ldots,n\}$. The polar invariant $\sigma_k(X,0)$ is defined by
$$\sigma_k(X,0) = \frac{1}{g_n^k} \int_{G_n^k} \sum_{i=1}^{N_P} \chi_i^P \cdot \Theta(K_i^p,0) dP.$$ We set $\sigma_0(X,0)=1$.}
\end{definition}

In \cite{ComteMerle} the authors defined another sequence of invariants attached to $X$, called the local Lipschitz-Killing invariants. 
\begin{definition}\label{Definition3.3}
{\rm Let $k \in \{0,\ldots,n\}$. The local Lipschitz-Killing invariant $\Lambda_k^{\rm loc}(X,0)$ is defined by
$$\Lambda_k^{\rm loc}(X,0) = \lim_{\epsilon \rightarrow 0} \frac{\Lambda_k(X \cap \mathbb{B}_\epsilon^n,X\cap \mathbb{B}_\epsilon^n)}{b_k \epsilon^k}.$$ }
\end{definition}
We note that $\Lambda_0^{\rm loc}(X,0)=1$. Then Comte and Merle proved linear kinematic formulas that relate the local Lipschitz-Killing invariants to the polar invariants.
\begin{theorem}[\cite{ComteMerle}, Theorem 3.1]\label{Theorem3.4}
For any  germ $(X,0) \subset (\mathbb{R}^n,0)$ of definable closed set, we have
$$ \left( \begin{array}{c}
\Lambda_1^{\rm loc}(X,0) \cr
\vdots \cr
\Lambda_n^{\rm loc}(X,0) \cr
\end{array} \right) =  \left( \begin{array}{cccc}
1 & m_1^2 & \ldots & m_1^n \cr
0 & 1 & \ldots & m_2^n \cr 
\vdots & \vdots & \ddots & \vdots \cr
0 & 0 & \ldots & 1 \cr
\end{array} \right) \cdot
\left( \begin{array}{c}
\sigma_1 (X,0) \cr
\vdots \cr
\sigma_n (X,0) \cr
\end{array} \right),$$
where $m_i^j = \frac{b_j}{b_{j-i} b_i} \binom{j}{i}  -\frac{b_{j-1}}{b_{j-1-i} b_i} \binom{j-1}{i} $, for $i+1 \le j \le n$.
\end{theorem}
If ${\rm dim} X=d$ then $\sigma_{d+1}(X,0)=\cdots= \sigma_n (X,0)=0$ and one recovers the local Cauchy-Crofton formula $\Lambda_d^{\rm loc}(X,0)= \sigma_d (X,0)$, previously proved by Comte \cite{Comte}. 

In \cite{DutertreGeoDedicata2012}, we also studied the asymptotic behavior of the Lipschitz-Killing measures. For $k=0,\ldots,n$, we considered the limits
$$\Lambda_k^{\rm lim}(X,0):= \lim_{\epsilon \rightarrow 0} \frac{\Lambda_k(X,X\cap \mathbb{B}_\epsilon^n)}{b_k \epsilon^k},$$ and we  showed the following theorem:
\begin{theorem}[\cite{DutertreGeoDedicata2012}, Theorem 5.1]\label{Theorem3.5}
For any  germ $(X,0) \subset (\mathbb{R}^n,0)$ of definable closed set, we have
$$\Lambda_0^{\rm lim}(X,0)=1-\frac{1}{2} \chi (\hbox{\em Lk}(X))-\frac{1}{2g_n^{n-1}} \int_{G_n^{n-1}} \chi (\hbox{\em Lk}(X \cap H)) dH.$$
Furthermore for $k \in \{1,\ldots,n-2 \}$, we have
$$\displaylines{
\qquad \Lambda_k^{\rm lim}(X,0)=- \frac{1}{2g_n^{n-k-1}} \int_{G_n^{n-k-1}} \chi(\hbox{\em Lk}(X \cap H)) dH \hfill \cr
\hfill + \frac{1}{2g_n^{n-k+1}} \int_{G_n^{n-k+1}} \chi(\hbox{\em Lk}(X \cap L)) dL, \qquad  \cr
}$$
and:
$$\Lambda_{n-1}^{\rm lim}(X,0)= \frac{1}{2g_n^2} \int_{G_n^2} \chi(\hbox{\em Lk}(X \cap H)) dH,$$
$$\Lambda_n^{\rm lim}(X,0) = \frac{1}{2g_n^1} \int_{G_n^1} \chi(\hbox{\em Lk}(X \cap H)) dH.$$
\end{theorem}

As a corollary, we obtained:
\begin{corollary}[\cite{DutertreGeoDedicata2012}, Corollary 5.2]\label{Corollary3.6}
For any  germ $(X,0) \subset (\mathbb{R}^n,0)$ of definable closed set, the equality
$1=\sum_{k=0}^n    \Lambda_k^{\rm lim}(X,0)$ holds.
\end{corollary}
We note that $\Lambda_k^{\rm lim} (X,0)$ differs from $\Lambda_k^{\rm loc}(X,0)$, because the link of $X$ does not have any contribution in the computation of $\Lambda_k^{\rm lim} (X,0)$.

In \cite{DutertreJofSingProcTrot}, we continued our study of the Lipschitz-Killing measures and established a local linear kinematic formula for $\Lambda_k^{\rm lim} (X,0)$, which was our inspiration for the present work. Let $H \in G_n^{n-k}$, $k \in \{1,\ldots, n \}$, and let $v$ be an element in $\mathbb{S}_{H^\perp}^{k-1}$. For $\delta >0$, we denote by $H_{v,\delta}$ the $(n-k)$-dimensional affine space $H+\delta v$ and we set 
$$\beta_0(H,v) = \lim_{\epsilon \rightarrow 0} \lim_{\delta \rightarrow 0} \Lambda_0(H_{\delta,v} \cap X, H_{\delta,v} \cap X \cap \mathbb{B}^n_\epsilon).$$
Then we set 
$$\beta_0 (H) =\frac{1}{s_{k-1}}\int_{\mathbb{S}_{H^\perp}^{k-1}} \beta_0 (H,v) dv.$$
\begin{theorem}[\cite{DutertreJofSingProcTrot}, Theorem 5.5]\label{Theorem3.7}
For $k \in \{1,\ldots,n \}$, we have 
$$ \Lambda_k^{\rm lim}(X,0) = \frac{1}{g_n^{n-k}} \int_{G_n^{n-k}} \beta_0 (H) dH.$$
\end{theorem}

We also proved local linear kinematic formulas that relate the limits $\Lambda_k^{\rm lim}(X,0)$ to the polar invariants. 
\begin{theorem}[\cite{DutertreJofSingProcTrot}, Theorem 5.6]\label{Theorem3.8}
For $k \in \{0,\ldots,n-1\}$, we have
$$\Lambda_k^{\rm lim}(X,0) = \sigma_k(X,0) -\sigma_{k+1}(X,0).$$
Furthermore, we have
$$\Lambda_n^{\rm lim}(X,0)= \sigma_n(X,0).$$
\end{theorem}

\section{Some preliminary topological and geometrical results}

Let $(X,0)$ and $(Y,0)$ be two germs of closed definable sets in $\mathbb{R}^n$. For convenience, we will work with two representatives of these germs that we denote by $X$ and $Y$ as well. We assume that these representatives $X$ and $Y$ are included in an open neighborhood $U$ of $0$. 

\subsection{A Gauss-Bonnet formula} Let $\{S_i\}_{i=0}^l$ be a Whitney stratification of $X$, where $S_0=\{0\}$ and $0 \in \overline{S_i}$ for $i \in \{1,\ldots,l\}$. Similarly let $\{T_j\}_{j=0}^m$ be a Whitney stratification of $Y$, where $T_0=\{0\}$ and $0 \in \overline{T_j}$ for $j= \{1,\ldots,m\}$. We assume each stratum to be connected. We introduce the following condition:
\begin{itemize}
\item Condition (1): for $i \in \{1,\ldots,l\}$ and for $j \in \{1,\ldots,m\}$, $S_i$ and $T_j$ intersect transversally (the case $S_i \cap T_j = \emptyset$ is possible).
\end{itemize}
If $X$ and $Y$ satisfy Condition (1) then $X \cap Y$ admits a Whitney stratification 
$X \cap Y = \sqcup_{k=0}^r R_k$, where $R_0=\{0\}$ and each $R_k$ is a connected component of an intersection $S_i \cap T_j$, $(i,j) \in \{1,\ldots,l\} \times \{1,\ldots,m\}$. 

Let $\widehat{X} \subset \mathbb{R}^{n+1}$ be the following definable set:
$$\widehat{X}= \left\{ (x,t) \in \mathbb{R}^{n+1} \ \vert \ x \in X \right\}.$$
It is included and closed in $U \times \mathbb{R}$. Let $v \in \mathbb{S}^{n-1}$ and let $\widehat{Y_v} \subset \mathbb{R}^{n+1}$ be the following definable set:
$$\widehat{Y_v} = \left\{ (y,t) \in \mathbb{R}^{n+1} \ \vert \
\exists y' \in Y \hbox{ such that } y=y'+tv \right\}.$$
It is included and closed in the open set 
$$\widehat{U_v} = \left\{ (u,t) \in \mathbb{R}^{n+1} \ \vert \
\exists u' \in U \hbox{ such that } u=u'+tv \right\}.$$
It is well-known (see \cite{GoreskyMacPherson}) that we can equip $\widehat{X}$ with a Whitney stratification in the following way: $\widehat{X}= \sqcup_{i=0}^l \widehat{S_i}$ where $\widehat{S_i}= \{ (x,t) \in \mathbb{R}^{n+1} \ \vert \ x \in S_i \}$. Similarly we can consider the following partition of $\widehat{Y_v}$: $\widehat{Y_v}=\sqcup_{j=0}^m \widehat{(T_j)_v}$ where 
$$\widehat{(T_j)_v}= \left\{ (y,t) \in \mathbb{R}^{n+1} \ \vert \
\exists y' \in T_j \hbox{ such that } y=y'+tv \right\}.$$
\begin{lemma}\label{Lemma4.1}
The partition $\widehat{Y_v}= \sqcup_{j=0}^m \widehat{(T_j)_v}$ gives a Whitney stratification of $\widehat{Y_v}$.
\end{lemma}
\proof With obvious notations, the partition $\widehat{Y}=\sqcup_{j=0}^m \widehat{T_j}$ induces a Whitney stratification of $\widehat{Y}$. Let $\phi : U \times \mathbb{R} \to \widehat{U_v}$ be defined by $\phi(u,t)=(u+tv,t)$. Then $\phi$ is a diffeomorphism, $\phi(\widehat{Y})= \widehat{Y_v}$ and $\phi ( \widehat{T_j})= \widehat{(T_j)_v}$ for $j \in \{0,\ldots,m\}$. This gives the result for Whitney's conditions are invariant by $C^1$-diffeomorphisms. \endproof
From now on, we will focus on the definable set $\widehat{X}\cap \widehat{Y_v}$. Let us denote if by $Z_v$. It is included and closed in the open set $(U \times \mathbb{R}) \cap \widehat{U_v}$. We introduce the following second condition:
\begin{itemize}
\item Condition (2): for $i \in \{0,\ldots,l\}$ and $j \in \{0,\ldots,m\}$, the strata $\widehat{S_i}$ and $\widehat{(T_j)_v}$ intersect transversally outside $(0,0)$.
\end{itemize}
If $v$ satisfies Condition (2) then $Z_v$ admits a Whitney stratification $Z_v = \sqcup_{l=0}^q Q_l$ where each $Q_l$ is a connected component of an intersection $\widehat{S_i} \cap \widehat{(T_j)_v}$. We note that necessarly $\widehat{S_0}\cap \widehat{(T_0)_v}=\{(0,0)\}$ and that we can put $Q_0=\{(0,0)\}$. 
\begin{lemma}\label{Lemma4.2}
Assume that $X$ and $Y$ satisfy Condition (1) and that $v$ satisfies Condition (2). Then the function 
$$\begin{array}{ccccc}
t_{\vert Z_v} & : & Z_v & \to & \mathbb{R} \cr
  & & (y,t) & \mapsto & t \cr
\end{array}$$
has an isolated stratified critical point at $(0,0)$.
\end{lemma}
\proof Let $Q$ be a stratum of $Z_v$ different from $\{(0,0)\}$. Since the critical points of $t_{\vert Q}$ lie in $\{t=0\}$, we can suppose that $Q$ is a connected component of $\widehat{S_i} \cap \widehat{(T_j)_v}$ with $i \not= 0$ and $j \not= 0$. Let us prove that $\{t=0\}$ intersects $\widehat{S_i} \cap \widehat{(T_j)_v}$ transversally. If it is not the case, then there is a point $p$ in $\widehat{S_i} \cap \widehat{(T_j)_v} \cap \{t=0\}$ such that $T_p ( \widehat{S_i} \cap \widehat{(T_j)_v} ) \subset \mathbb{R}^n$. But it is not difficult to check that $\{t=0\}$ intersects $\widehat{S_i}$ and $\widehat{(T_j)_v}$ transversally, so $T_p S_i = T_p \widehat{S_i} \cap \mathbb{R}^n $ and $T_p T_j=T_p \widehat{(T_j)_v} \cap \mathbb{R}^n$. Moreover, $S_i$ and $T_j$ intersect transversally and so $T_p(S_i \cap T_j)= T_p S_i \cap T_p T_j$. Similarly, $T_p (\widehat{S_i} \cap \widehat{(T_j)_v})= T_p \widehat{S_i} \cap T_p \widehat{(T_j)_v}$. We get that $T_p (S_i \cap T_j)=T_p (\widehat{S_i} \cap \widehat{(T_j)_v})$. This is not possible, for ${\rm dim} \widehat{S_i}={\rm dim} S_i +1$, ${\rm dim} \widehat{(T_j)_v}={\rm dim} T_j +1$ and these two intersections are transverse in $\mathbb{R}^n$ and $\mathbb{R}^{n+1}$. \endproof
We can apply Corollary 4.6 in \cite{DutertreJofSingProcTrot} to $t$ and $Z_v$.
\begin{corollary}\label{Corollary4.3}
Assume that $X$ and $Y$ satisfy Condition (1) and that $v$ satisfies Condition (2). Then we have
$$\displaylines{
\quad \lim_{\epsilon \to 0} \lim_{\delta \to \pm 0} 
\Lambda_0 \left(Z_v \cap \{t=\delta \}, Z_v \cap \{t=\delta \} \cap \mathbb{B}_\epsilon^{n+1} \right) \hfill \cr
\qquad \qquad= \lim_{\epsilon \to 0} \lim_{\delta \to \pm 0}
\chi \left(Z_v \cap \{t=\delta \} \cap \mathbb{B}_\epsilon^{n+1} \right) \hfill \cr
\hfill -\frac{1}{2} \chi \left( {\rm Lk} (X\cap Y) \right)
-\frac{1}{2 s_{n-1}}\int_{\mathbb{S}^{n-1}} \chi \left({\rm Lk} (X \cap Y \cap \{u^*=0 \} \right)du. \quad \cr 
}$$
\end{corollary}
\proof We just have to show that 
$$\frac{1}{2 s_{n-1}}\int_{\mathbb{S}^{n-1}} \chi \left({\rm Lk} (X \cap Y \cap \{u^*=0 \} \right)du =\frac{1}{2 s_n}\int_{\mathbb{S}^n} \chi \left({\rm Lk} (X \cap Y \cap \{u^*=0 \} \right)du.$$
But since $X \cap Y$ is included in $\mathbb{R}^n$, the method given in the proof of Corollary 5.1 \cite{DutertreAdvGeo2008} applies here. \endproof

Let us go back now to the sets $X$ and $Y$. We denote the definable set $X \cap (Y +\delta v)$ by $Z_{v,\delta}$.
\begin{lemma}\label{Lemma4.4} There exists $\epsilon_0 >0$ such that for $0< \epsilon \le \epsilon_0$, there exists $\delta_\epsilon >0$ such that for $0< \delta \le \delta_\epsilon$, the topological type of $Z_{v,\delta} \cap \mathbb{B}_\epsilon^n$ does not depend on the choice of the couple $(\epsilon,\delta)$. Moreover for $0< \delta \ll \epsilon \ll 1$, 
$ Z_{v,\delta} \cap \mathbb{B}_\epsilon^n$ and
$ Z_v \cap \{t=\delta\} \cap \mathbb{B}_\epsilon^{n+1} $ are homeomorphic.
\end{lemma}
\proof Let $\rho(x,t)=\sqrt{x_1^2+\cdots+x_n^2}+t$. Then $Z_{v,\delta} \cap \mathbb{B}_\epsilon^n$ is homeomorphic to $Z_v \cap \{t=\delta\} \cap \{\rho \le \epsilon + \delta \}$. Let $\pi : Z_v \cap \{t \ge 0 \} \to \mathbb{R}$ be the mapping defined by $\pi(x,t)= (\rho(x,t),t)$ and let $\Delta \subset \mathbb{R}^2$ be its (stratified) discriminant. It is a definable curve in a neighborhood of $0 \in \Delta$ in $\mathbb{R}^+ \times \mathbb{R}^+$. The following function 
$$\begin{array}{ccccc}
r & : & \mathbb{R}^+ & \to & \mathbb{R}^+ \cr
 & & \nu & \mapsto & {\rm inf} \{t \ \vert \ (\nu,t) \in \Delta \} \cr
\end{array}$$
is a definable function defined in a neighborhood of $0$. Note that $r(0)=0$ and $r(\nu)>0$ for $\nu>0$ close to $0$. By the Monotonicity Theorem (see Theorem 2.1 in \cite{CosteOmin} or \cite{DriesMiller}, 4.1), we can assume that $r$ is continuous and  increasing on a small interval $]0,\nu_0[$. Let $(\epsilon,\delta)$ be a couple such that $0< \epsilon < \nu_0$, $0< \delta < r(\epsilon)$ and $Z_v \cap \{t=\delta\} \cap \{ \rho \le \epsilon \}$ has the topological type of the Milnor fibre of $t_{\vert Z_v}$. By taking $\delta$ smaller if necessary, we can assume that $\epsilon + \delta < \nu_0$. Since $r$ is increasing, $\delta < r(\epsilon) < r(\epsilon+\delta)$ and so $Z_v \cap \{t=\delta \} \cap \{\rho \le \epsilon+\delta \}$ is homeomorphic to $Z_v \cap \{t=\delta \} \cap \{\rho \le \epsilon \}$. We conclude with the results of Section 3. \endproof
A similar result is true for negative values of $t$ replacing $\sqrt{x_1^2+\cdots+x_n^2}+t$ with $\sqrt{x_1^2+\cdots+x_n^2}-t$. We can state the infinitesimal Gauss-Bonnet formula for $Z_{v,\delta}$. 
\begin{lemma}\label{Lemma4.5}
Assume that $X$ and $Y$ satisfy Condition (1) and that $v$ satisfies Condition (2). Then we have
$$\displaylines{ \qquad \lim_{\epsilon \to 0} \lim_{\delta \to \pm 0} 
\Lambda_0 \left(Z_{v,\delta}, Z_{v,\delta} \cap \mathbb{B}_\epsilon^{n} \right) 
 = 
\lim_{\epsilon \to 0} \lim_{\delta \to \pm 0} \chi \left(Z_{v,\delta} \cap \mathbb{B}_\epsilon^{n} \right) \hfill \cr
\hfill -\frac{1}{2} \chi \left( {\rm Lk} (X\cap Y) \right)
-\frac{1}{2 s_{n-1}}\int_{\mathbb{S}^{n-1}} \chi \left({\rm Lk} (X \cap Y \cap \{u^*=0 \} )\right)du. \qquad \cr
}$$
\end{lemma}
\proof Let $i : \mathbb{R}^n \to \mathbb{R}^{n+1}$, $x \mapsto (x,\delta)$. Since $i$ is a definable isometry, by Theorem 5.0 in \cite{Fu94} or Proposition 9.2 in \cite{BroeckerKuppe}, we have 
$$\Lambda_0^n \left(Z_{v,\delta}, Z_{v,\delta} \cap \mathbb{B}^n_\epsilon \right)=
\Lambda_0^{n+1}  \left(Z_v \cap \{t=\delta \} , Z_v \cap \{t=\delta \} \cap \{\rho \le \epsilon+\delta \} \right).$$
Here we suppose that $\delta >0$, $\Lambda_0^n$ (resp. $\Lambda_0^{n+1}$) stands for the Gauss-Bonnet measure in $\mathbb{R}^n$ (resp. $\mathbb{R}^{n+1}$) and $\rho(x,t)=\sqrt{x_1^2+\cdots+x_n^2}+t$ 

As explained in Section 3, for $u$ generic in $\mathbb{S}^n$ and for $\epsilon >0$ sufficiently small, there exists $\delta_{\epsilon,u}$ such that for $0< \delta \le \delta_{\epsilon,u}$, the critical points of $u^*$ and $-u^*$ in $Z_{v,\delta} \cap \{t=\delta\}$ actually lie in $Z_{v,\delta} \cap \{t=\delta\} \cap \{ \rho \le \frac{\epsilon}{4} \}$, hence there are not in $Z_{v,\delta} \cap \{t=\delta\} \cap \{\epsilon \le \rho \le \epsilon+\delta \}$. Thanks to this observation, we can conclude that 
$$\lim_{\epsilon \to 0} \lim_{\delta \to 0} 
\Lambda_0^{n+1} \left(Z_v \cap \{t=\delta \}, Z_v \cap \{t=\delta \} \cap \{\rho \le \epsilon+\delta\} \right)$$ $$= \lim_{\epsilon \to 0} \lim_{\delta \to 0} 
\Lambda_0^{n+1} \left(Z_v \cap \{t=\delta \}, Z_v \cap \{t=\delta \} \cap \{\rho \le \epsilon\} \right).$$
It is enough to apply Corollary \ref{Corollary4.3} and the comments of Section 3 on the choice of the distance function to the origin to get the result. \endproof

\subsection{A useful lemma}
We continue this section with a remark. Instead of translating $Y$, we can translate $X$, intersect this translated set with $Y$ and obtain another Milnor fibre $Y \cap (X+\delta v) \cap \mathbb{B}_\epsilon^n$, $0 < \vert \delta \vert \ll \epsilon \ll 1$.
\begin{lemma}\label{Lemma4.6}
There exists $\epsilon_0 >0$ such that for $0< \epsilon \le \epsilon_0$, there exists $\delta_\epsilon >0$ such that $0< \delta \le \delta_\epsilon$, $X \cap (Y + \delta v) \cap \mathbb{B}_\epsilon^n$ and $Y\cap (X-\delta v)\cap \mathbb{B}_\epsilon^n$ are homeomorphic.
\end{lemma}
\proof Let $\rho_v(x,t)=\sqrt{(x_1-tv_1)^2+\cdots+(x_n-tv_n)^2}+t$. By the results of Section 3, we know that there exists $\epsilon_0>0$ such that for $0< \epsilon \le \epsilon_0$, there exists $\delta_\epsilon >0$ such that for $0< \delta \le \delta_\epsilon$, the topological type of $Z_v \cap \{t=\delta\} \cap \{\rho_v \le \epsilon \}$ does not depend on the couple $(\epsilon,\delta)$ and is the topological type of the positive Milnor fibre of $t_{\vert Z_v}$. On the other hand, the set $X \cap (Y+\delta v) \cap (\mathbb{B}_\epsilon^n+\delta v)$ is homeomorphic to $Z_v \cap \{t=\delta \} \cap \{\rho_v \le \epsilon + \delta\}$. The same method as the one used in Lemma \ref{Lemma4.4} shows that for $0< \epsilon \le \epsilon_0$ and $0< \delta \le \delta_\epsilon$ small enough, 
$Z_v \cap \{t=\delta \} \cap \{\rho_v \le \epsilon + \delta\}$ is homeomorphic to $Z_v \cap \{t=\delta \} \cap \{\rho_v \le \epsilon\}$. We conclude that $Z_{v,\delta} \cap \mathbb{B}_\epsilon^n$ is homeomorphic to $Z_{v,\delta} \cap (\mathbb{B}_\epsilon^n +\delta v)$ for $0< \delta \ll \epsilon \ll 1$. But $$\Big[ X \cap (Y+\delta v) \cap (\mathbb{B}_\epsilon^n + \delta v) \Big]-\delta v =(X-\delta v) \cap Y \cap \mathbb{B}_\epsilon^n.$$ \endproof

\subsection{Genericity of Conditions (1) and (2)}
We prove the genericity of Conditions (1) and (2). To prove the genericity of Condition (1), we need some auxiliary lemmas.
\begin{lemma}\label{Lemma4.7}
Let $x \in \mathbb{R}^n$ be a non-zero vector. We have 
$$\left\{ H x \ \vert \ H \in M_n (\mathbb{R}) \hbox{ such that }
{}^t \! H=-H \right\}= x^\perp.$$
\end{lemma}
\proof It is clear that if $H$ is an antisymmetric matrix, then $Hx \in x^\perp$. Let us write $x=(x_1,\ldots,x_n)$ and let $a=(a_1,\ldots,a_n) \in x^\perp$. Since $x \not= (0,\ldots,0)$ then there exists $k$ such that $x_k \not= 0$. Then we can construct $H=(h_{ij})$ in the following way:
$$h_{kj}=-\frac{a_j}{x_k}, \ h_{jk}=-h_{kj} \hbox{ for } j \not= k,$$
and putting $h_{ij}=0$ for the other coefficients. Then $H$ is antisymmetric and $Hx=a$. \endproof

\begin{lemma}\label{Lemma4.8}
Let $f : \mathbb{R}^n \to \mathbb{R}^k$ with $1 \le k \le n-1$ be a $C^1$ mapping and let $F$ be the mapping defined by
$$\begin{array}{ccccc}
F & : & M_n(\mathbb{R}) & \to & S_n (\mathbb{R}) \times \mathbb{R}^k \cr
 &  & A & \mapsto & ({}^t \!AA, f(Ax)) , \cr
\end{array} $$
where $x$ is a non-zero vector. If $A \in SO(n)$ and $Df(Ax)_{\vert (Ax)^\perp} : (Ax)^\perp \to \mathbb{R}^k$ is surjective then $DF(A)$ is a surjection.
\end{lemma}
\proof We have $DF(A)(H)=({}^t \!AH+{}^t \!HA, Df(Ax)(Hx))$. Let $(Y,\alpha) \in S_n(\mathbb{R}) \times \mathbb{R}^k$, we have that 
$${}^t \!A \left( \frac{1}{2}AY \right) +
{}^t \!  \left( \frac{1}{2}AY \right) A =Y.$$ Let $\beta= Df(Ax) (\frac{1}{2}AY x)$. We have to find $H$ such that ${}^t \!A H + {}^t \!H A=0$ and $Df(Ax)(Hx)=\alpha-\beta$. Since $Df(Ax)_{\vert (Ax)^\perp}$ is a surjection, by the previous lemma, there exists an antisymmetric matrix $L$ such that $Df(Ax)(LAx)=\alpha -\beta$. We take $H=LA$. \endproof

\begin{lemma}\label{Lemma4.9}
Let $T \subset \mathbb{R}^n$ be a $C^1$ definable submanifold of dimension $d$ such that $0 \in \overline{T}$. Then there exists a neighborhood $U_T$ of $0$ such that for $x \in T \cap (U_T \setminus \{0\} )$, ${\rm dim} (T_x T \cap x^\perp) \le d-1$.
\end{lemma}
\proof If it is not the case, then there exists an injective $C^1$ definable map $\mu : [0,\nu) \to \overline{T}$ such that $\mu(0)=0$ and for $t \not= 0$, $T_{\mu(t)} T \subset \mu (t)^\perp$, hence $\mu (t) \perp T_{\mu (t)} T$. Since $\mu'(t) \in T_{\mu (t)} T$, we get that $\langle \mu (t), \mu ' (t) \rangle=0$. This implies that $\langle \mu (t) ,\mu (t) \rangle' =0$ and that $\vert \mu (t) \vert =0$, which is not possible. \endproof

\begin{lemma}\label{Lemma4.10}
There exists a definable subset $\Sigma_{X,Y} \subset SO(n)$ of positive dimension such that for $\gamma \notin \Sigma_{X,Y}$, $X$ and $\gamma Y$ satisfy Condition (1).
\end{lemma}
\proof We note first that $\{ \gamma T_j \}_{j=0}^m$ is a Whitney stratification of $\gamma Y$. Let $S\not= \{0\}$ be a stratum of $X$ and let $T \not= \{0\}$ be a stratum of $Y$. We have to prove that there exists a definable subset $\Sigma_{S,T} \subset SO(n)$ of positive codimension such that $S$ and $\gamma T$ intersect transversally in a neighborhood of $0$. If ${\rm dim} S =n$ or ${\rm dim} T=n$, there is nothing to prove so we can assume that $e:={\rm dim} S \le n-1$ and $d:={\rm dim} T \le n-1$. Let $M$ be the following definable set:
$$ M =\left\{ (p,\gamma) \in (U'\setminus \{0\}) \times O(n) \ \vert \ 
p \in S \cap \gamma^{-1} T \right\},$$
where $U'$ is an open definable neighborhood of $0$, included in $U \cap U_T$ and such that $\gamma U' \subset U'$ for all $\gamma \in O(n)$. Let us prove that $M$ is a definable submanifold of $\mathbb{R}^n \times M_n (\mathbb{R})$. 

Let $(p,\gamma)$ be a point in $M$. There is an open neighborhood $\mathcal{V}$ of $(p,\gamma)$ in $\mathbb{R}^n \times M_n(\mathbb{R})$ such that in $\mathcal{V}$, $M$ is defined $G(x,A)=(0,I_n,0)$, where 
$$\begin{array}{ccccc}
G & : & \mathbb{R}^n \times M_n(\mathbb{R}) & \to & \mathbb{R}^{n-e} \times S_n(\mathbb{R}) \times \mathbb{R}^{n-d} \cr
 &  & (x,A) & \mapsto & (g(x),{}^t \!AA, f(Ax)) , \cr
\end{array}$$
and where $g$ is a definable submersion such that $S$ is locally defined by $g(x)=0$ in a neighborhood of $p$ and $f$ is a definable submersion such that $T$ is locally defined by $f(x)=0$ in a neighborhood of $\gamma p$. 

Since $p$ belongs to $U'$, $\gamma p \in U_T \cap T$ and so ${\rm dim} (T_{\gamma p} T \cap (\gamma p)^\perp ) \le d-1$, which implies that $Df(\gamma p)_{\vert (\gamma p)^\perp} : (\gamma p)^\perp \to \mathbb{R}^{n-d}$ is a surjection. By Lemma \ref{Lemma4.8}, the mapping
$$\begin{array}{ccccc}
F & : & M_n (\mathbb{R}) & \to & S_n(\mathbb{R}) \times \mathbb{R}^{n-d} \cr
 &  & A & \mapsto & ({}^t \!AA,f(Ap) )\cr
\end{array}$$
is a submersion at $\gamma$. Therefore the submatrix of the jacobian matrix of $G$ at $(p,\gamma)$ formed by the partial derivatives of $F$ with respect $A$ has maximal rank. But the submatrix of the jacobian matrix of $G$ at $(p,\gamma)$ formed by the partial derivatives of $g$ with respect to $x$ has also maximal rank. We conclude that $G$ is a submersion at $(p,\gamma)$ and that $M$ is a definable submanifold of dimension 
$$n+n^2-2n-\frac{n(n+1)}{2}+(e+d)= \frac{n(n-1)}{2} + (e+d)-n.$$
Let $\pi : M \to O(n)$ be the natural projection. By Sard's theorem (see \cite{BochnakCosteRoy}), its discriminant $\Delta$ is a definable subset of positive codimension. Let $\mathcal{T} : O(n) \to O(n)$ be the definable diffeomorphism given by $\mathcal{T} (A)={}^t \!A$. It is enough to take $\Sigma_{S,T} =\mathcal{T} (\Delta) \cap SO(n)$.
\endproof

\begin{lemma}\label{Lemma4.11}
Assume that $X$ and $Y$ satisfy Condition (1). There exists a definable subset $\Gamma_{X,Y} \subset \mathbb{S}^{n-1}$ of positive codimension such that for $v \notin \Gamma_{X,Y}$, $v$ satisfies Condition (2).
\end{lemma}
\proof Let $S$ be a stratum of $X$ and let $T$ be a stratum of $Y$. We have to prove that there exists a definable subset $\Gamma_{S,T} \subset \mathbb{S}^{n-1}$ of positive codimension such that $\widehat{S}$ and $\widehat{T_v}$ intersect transversally outside $(0,0)$, in a neighborhood of $(0,0)$. 

If $S=\{0\}$ and $T=\{0\}$, then $\widehat{S} \cap \widehat{T_v}= \{(0,0)\}$ and there is nothing to prove. If ${\rm dim} S =n$ or ${\rm dim} T=n$, there is nothing to prove neither. Let us treat first the case $0< e:={\rm dim} S \le n-1$ and $0< d:={\rm dim} T \le n-1$. Let $M$ be the following definable set:
$$ M= \left\{ (p,\tau,\nu) \in U \times \mathbb{R} \times (\mathbb{R}^n \setminus \{0\}) 
\ \vert \ (p,\tau) \in \widehat{S} \cap \widehat{T_\nu} \right\}.$$
Let us prove that $M$ is a definable submanifold. Let $(p,\tau,\nu)$ be a point in $M$. There is an open neighborhood $\mathcal{V}$ of $(p,\tau,\nu) \in U \times \mathbb{R} \times (\mathbb{R}^n \setminus \{0\}) $ such that in $\mathcal{V}$, $M$ is defined by $g(x)=0$ and $f(x-tv)=0$, where $g$ is a definable submersion such that $S$ is locally defined by $g(x)=0$ in a neighborhood of $p$ and $f$ is a definable submersion such that $T$ is locally defined by $f(x)=0$ in a neighborhood of $p-\tau \nu$. It is easy to check that the Jacobian matrix of the mapping $(g,f)$ has maximal rank at $(p,\tau,\nu)$ if $\tau \not= 0$. If $\tau =0$ this is also the case by Condition (1). Therefore $M$ is a definable submanifold of dimension $2n+1-2n+(e+d)=(e+d)+1$. Let $\pi : M \to \mathbb{R}^ n$ be the projection $\pi(p,\tau,\nu)=\nu$. By Sard's theorem (see \cite{BochnakCosteRoy}), its discriminant $\Delta$ is a definable subset of positive codimension. We take $\Gamma_{S,T}=\mathbb{S}^{n-1} \setminus \Delta$. The remaining two cases are proved with the same method.
\endproof 

\subsection{The link of $Z_v \cap \{ t\ge 0\}$}
We study the link of the set $Z_v \cap \{ t\ge 0\}$. We still assume that $X$ and $Y$ satisfy Condition (1) and that $v$ satisfies Condition (2). For $(x,t) \in \mathbb{R}^{n+1}$, we set $\omega(x,t)=\vert x \vert$. Let $\epsilon_0>0$ be such that $\mathbb{B}^n_{\epsilon_0} \subset U$. We set 
$$\Gamma_Y =\left\{ (\frac{1}{u}y,u) \ \vert \ y \in Y, u \in ]0,\epsilon_0[ \right\} \subset \mathbb{R}^{n+1}.$$
We recall that the tangent cone of $Y$at $0$ is $C_0 Y = \overline{\Gamma_Y} \cap (\mathbb{R}^n \times \{0\})$. It is a closed and conic definable set in $\mathbb{R}^n$. 

\begin{lemma}\label{Lemma4.12}
If $v \notin \left(-C_0(Y) \right) \cap \mathbb{S}^{n-1}$ then there exist $\epsilon_v >0$  and $a>0$ such that the inclusion
$$\widehat{Y_v} \cap \{ t \ge 0 \} \cap \mathbb{B}_{\epsilon_v}^{n+1} \subset
\left\{ (x,t) \in \mathbb{R}^{n+1} \ \vert \ \vert x \vert \ge at \right\} \cap \mathbb{B}_{\epsilon_v}^{n+1}$$ holds.
\end{lemma}
\proof If it is not the case then we can find  a sequence of points $(x_n,t_n)_{n \in \mathbb{N}}$ in $\widehat{Y_v}\setminus \{ t=0\}$ such that $(x_n,t_n) \to (0,0)$ and $\lim_{n \to +\infty} \frac{\vert x_n \vert}{t_n}=0$. We have 
$$\vert x_n-t_nv \vert^2 =\vert x_n \vert^2+t_n^2-2t_n \langle x_n ,v \rangle,$$
and so
$$\frac{\vert x_n-t_nv \vert^2}{t_n^2}=\frac{\vert x_n \vert^2}{t_n^2}+1-2 \langle \frac{x_n}{t_n},v \rangle.$$
Since $\vert \langle \frac{x_n}{t_n},v \rangle \vert \le \frac{\vert x_n \vert}{t_n}$, we find that $\lim_{n\to +\infty} \frac{\vert x_n-t_nv \vert}{t_n}=1$. Since $\frac{\vert x_n \vert}{\vert x_n-t_nv\vert}=\frac{\vert x_n \vert}{t_n} \times \frac{t_n}{\vert x_n-t_nv \vert}$, we find that $\lim_{n \to +\infty} \frac{\vert x_n \vert}{\vert x_n-t_n v \vert}=0$. Therefore we see that $$\lim_{n \to +\infty} \frac{x_n-t_nv}{\vert x_n-t_n v \vert}=-v,$$ which implies that $-v$ belongs to $C_0(Y) \cap \mathbb{S}^{n-1}$. \endproof

\begin{corollary}\label{Corollary4.13}
If $v \notin \left(-C_0(Y) \right) \cap \mathbb{S}^{n-1}$ then for $0< \epsilon \ll 1$, the link of $\widehat{Z_v} \cap \{t \ge 0\}$ is homemorphic to $\widehat{Z_v} \cap \{t \ge 0 \}
\cap \left\{ (x,t) \ \vert \ \vert x \vert =\epsilon \right\}$.
\end{corollary}
\proof It is enough to show the implication 
$$\left\{ \begin{array}{c}
(y_0,t_0) \in \widehat{Y_v} \cap \{t \ge 0 \} \cap \mathbb{B}_{\epsilon_v}^{n+1} \cr
\omega (y_0,t_0)=0 \cr
\end{array} \right. \Rightarrow (y_0,t_0)=(0,0). $$
But if $\omega (y_0,t_0)=0$ then $\vert y_0 \vert=0$ and, by the previous lemma, $\vert y_0 \vert \ge a t_0$ for some $a >0$. \endproof

If ${\rm dim} C_0 (Y) \le n-1$ then the previous corollary holds for almost all $v$ in $\mathbb{S}^{n-1}$ since ${\rm dim} C_0 (Y) \cap \mathbb{S}^{n-1} < n-1$. 

In the rest of this subsection, we assume that ${\rm dim} C_0(Y)=n$. This implies that ${\rm dim} Y=n$. Let us denote by $Y'$ the union of the strata of $Y$ of dimension less than or equal to $n-1$. It is a closed definable set of dimension less than or equal to $n-1$, if not empty. 
We need auxiliary lemmas.
\begin{lemma}\label{Lemma4.14}
Let $S \subset \mathbb{R}^n$ be a definable open subset such that $0 \in \overline{S}$ and that ${\rm dim} C_0 S=n$. If $v$ is in $C_0 S \setminus C_0(\overline{S} \setminus S)$ then there is $\alpha >0$ such that $]0,\alpha]\cdot v \subset S$.
\end{lemma}
\proof We note that $C_0 S \setminus C_0(\overline{S} \setminus S)$ is not empty because ${\rm dim} C_0(S)=n$. Let $v \in C_0 S \setminus C_0(\overline{S} \setminus S)$ (note that necessarly $v \not= 0$). Let us suppose that for all $\alpha >0$, $]0,\alpha] \cdot v$ is not included in $S$. Hence we can construct a sequence $(z_n)_{n\in \mathbb{N}}$ in the complement ${}^c \! S$ of $S$ such that $(z_n)$ tends to $0$ and $\frac{z_n}{\vert z_n \vert}=v$. This implies that $v$ belongs to $C_0({}^c \! S)$. Since $v$ belongs to $C_0(S)$ as well, there exists a sequence of points $(x_n)_{n \in \mathbb{N}}$ in $S$ such that $\frac{x_n}{\vert x_n \vert}$ tends to $v$. 

Let $a \in S$ and let $b \in {}^c \! S$. Let $[a,b]$ be the segment with extremities $a$ and $b$, i.e.
$$[a,b]= \left\{ z \ \vert \ z=\lambda a +(1-\lambda) b, \ \lambda \in [0,1] \right\}.$$
Since $S$ is open, there exists $0 < \eta \le 1$ such that $[a,a+\eta (b-a)[$ is included in $S$. Let $\eta_a$ be the supremum of such $\eta$'s. The point $a+\eta_a (b-a)$ lies in $\overline{S}\setminus S$. Actually if $a+\eta_a (b-a)$ is in $S$ then there exists $\eta' > \eta_a$ such that $a+\eta' (b-a) \in S$, because $S$ is open. We conclude that the segment $[a,b]$ intersects $\overline{S}\setminus S$. 

For each $n \in \mathbb{N}$, let $y_n$ be a point in $[x_n,z_n] \cap (\overline{S} \setminus S)$. If there is a subsequence $y_{\tau(n)}$ such that $y_{\tau (n)}= 0$, then $\frac{x_{\tau (n)}}{\vert x_{\tau (n)} \vert}=-v$, which is not possible for $\frac{x_{\tau (n)}}{\vert x_{\tau (n)} \vert}$ tends to $v$. Therefore we can assume that $y_n \not= 0$ for $n \in \mathbb{N}$, and write
$$\frac{y_n}{\vert y_n \vert}= \frac{\lambda_n \vert z_n \vert}{\vert y_n \vert}\cdot v+
\frac{(1-\lambda_n) \vert x_n \vert}{\vert y_n \vert}\cdot \frac{x_n}{\vert x_n \vert},$$
where $\lambda_n \in [0,1]$. For simplicity we rewrite this equality in the following way:
$$\frac{y_n}{\vert y_n \vert}= \alpha_n \cdot v+\beta_n \cdot \frac{x_n}{\vert x_n \vert},$$
with $\alpha_n,\beta_n \ge 0$. Since $\langle v, \frac{x_n}{\vert x_n \vert} \rangle \to 1$, there is $n_0$ such that for $n \ge n_0$, $\langle v, \frac{x_n}{\vert x_n \vert} \rangle \ge \frac{1}{2}$. This implies that for $n \ge n_0$, $0 \le \alpha_n^2+\beta_n^2+\alpha_n \beta_n \le 1$ and so the sequences $(\alpha_n)_{n \in \mathbb{N}}$ and $(\beta_n)_{n \in \mathbb{N}}$ are bounded. Therefore, taking a subsequence if necessary, we can assume that $\alpha_n $ tends to $\alpha \ge 0$ and $\beta_n$ tends to $\beta \ge 0$. Hence $\frac{y_n}{\vert y_n \vert}$ tends to $(\alpha+\beta)v$, where $\alpha+\beta=1$ for the limit of $\frac{y_n}{\vert y_n \vert}$ is a unit vector. We see that $v$ belongs to $C_0(\overline{S}\setminus S)$, which is not possible by hypothesis. We conclude that there is $\alpha >0$ such that $]0,\alpha] \cdot v \subset S$. \endproof
\begin{lemma}\label{Lemma4.15}
Let $W \subset \mathbb{R}^n$ be a closed definable set equipped with a Whitney stratification. Suppose that $0 \in W$ and that $0$ lies in a stratum of dimension greater than or equal to 1. Let $g : W \to \mathbb{R}$ be a definable function, restriction of a $C^2$ definable function, such that $0$ is not a stratified critical point of $g$. 
Let $f : W \cap \{ g \le 0\} \to \mathbb{R}$ be a definable function, restriction of a $C^2$ definable function, such that $f(0)=0$ and $0$ is a local strict maximum of $f$. Then ${\rm ind}(f,W \cap \{g \le 0\},0)=0$.
\end{lemma} 
\proof By Lemma 3.3 in \cite{DutertreJofSingNonIsol}, $\chi \left( {\rm Lk} (W \cap \{ g \le 0 \} ) \right)=1$. 
Since $-f_{\vert W \cap \{ g \le 0 \}  }$ is a distance function to the origin, we can write 
$${\rm ind}(f,W \cap \{g \le 0\},0) =1-\chi \left( W \cap \{g \le 0 \} \cap \{ -f \le \epsilon\} \cap 
\{ f=-\delta \} \right),$$
with $0< \delta <\epsilon \ll 1$. But $W \cap \{g \le 0 \} \cap \{ -f \le \epsilon\} \cap 
\{ f=-\delta \}$ is homeomorphic to the link of $W \cap \{ g \le 0 \}$ at $0$. \endproof

Let us choose $v$ in $\left( -C_0(Y) \right) \cap \mathbb{S}^{n-1}$ such that $v \notin \left(-C_0(Y') \right) \cap \mathbb{S}^{n-1}$. By Lemma \ref{Lemma4.12}, there exists $\epsilon_v >0$ and $a >0$ such that 
$$\widehat{Y'_v} \cap \{t \ge 0 \} \cap \mathbb{B}_{\epsilon_v}^{n+1} \subset 
\left\{ (x,t) \ \vert \ \vert x \vert \ge a t \right\} \cap \mathbb{B}_{\epsilon_v}^{n+1}.$$
\begin{lemma}\label{Lemma4.16}
Under these assumptions, we have
$$\chi  \left( X \cap (Y+\frac{2\epsilon}{a}v ) \cap \mathbb{B}_\epsilon^n \right)= 1,$$
for $0< \epsilon \ll 1$.
\end{lemma}
\proof For $0 < \epsilon \ll 1$, we set $r_\epsilon =\sqrt{\epsilon^2 + \left( \frac{2\epsilon}{a}\right)^2}$. The set $X \cap (Y+\frac{2\epsilon}{a}) \cap \mathbb{B}_\epsilon^n$ is the intersection $Z_v \cap \{ t =\frac{2\epsilon}{a} \} \cap \mathbb{B}_{r_\epsilon}^{n+1}$. If $(x_0,t_0)$ lies in $Z_v \cap \{ t \ge \frac{2\epsilon}{a} \} \cap \mathbb{B}_{r_\epsilon}^{n+1}$ then $\vert x_0 \vert ^2  \le \epsilon^2 + \left( \frac{2 \epsilon}{a} \right)^2-t_0^2$ and $\frac{\vert x_0 \vert^2}{t_0^2}\le \frac{\epsilon^2}{t_0^2} \le \left( \frac{a}{2} \right)^2$. Therefore if $\epsilon$ is sufficiently small, $Z_v \cap \{ t \ge \frac{2\epsilon}{a} \} \cap \mathbb{B}_{r_\epsilon}^{n+1} \cap \widehat{Y'_v}= \emptyset$. 

The function $t_{\vert Z_v}$ has an isolated  stratified critical value at $0$. If $\epsilon$ is small enough, the stratified critical points of $t_{\vert Z_v \cap \{ t > \frac{2\epsilon}{a}\} \cap \mathbb{B}_{r_\epsilon}^{n+1}}$ lie in $Z_v \cap \{ t > \frac{2\epsilon}{a} \}\cap \mathbb{S}_{r_\epsilon}^n$. But $\widehat{Y'_v}$ does not intersect $Z_v \cap \{ t > \frac{2\epsilon}{a} \}\cap \mathbb{B}_{r_\epsilon}^{n+1}$, so $Z_v \cap \{t > \frac{2\epsilon}{a} \} \cap \mathring{\mathbb{B}_{r_\epsilon}^{n+1}}$ (resp. $Z_v \cap \{t > \frac{2\epsilon}{a} \} \cap \mathbb{S}_{r_\epsilon}^{n}$) is stratified by strata of the form $\widehat{S} \cap \widehat{T_v} \cap \{t > \frac{2\epsilon}{a} \} \cap \mathring{\mathbb{B}_{r_\epsilon}^{n+1}}$ (resp. $\widehat{S} \cap \widehat{T_v} \cap \{t > \frac{2\epsilon}{a} \} \cap \mathbb{S}_{r_\epsilon}^{n}$), where $S$ is a stratum of $X$ and $T$ is a stratum of $Y$ of dimension $n$. This means that $Z_v \cap \{ t > \frac{2\epsilon}{a}\} \cap \mathring{\mathbb{B}_{r_\epsilon}^{n+1}}$ (resp. $Z_v \cap \{t > \frac{2\epsilon}{a} \} \cap \mathbb{S}_{r_\epsilon}^{n}$) is stratified by open subsets of strata of the form $\widehat{S} \cap \{ t  > \frac{2\epsilon}{a} \} \cap \mathring{\mathbb{B}_{r_\epsilon}^{n+1}}$ (resp. $\widehat{S} \cap \{ t  > \frac{2\epsilon}{a} \} \cap \mathbb{S}_{r_\epsilon}^{n}$). 

Such a stratum $\widehat{S}$ is a product $S \times ]-\epsilon',\epsilon"[$ where $\epsilon',\epsilon" >0$ are small and $S$ is a stratum of $X$. A point $(x_0,t_0)$ is a critical point of $t_{\vert \widehat{S} \cap \mathbb{S}_{r_\epsilon}^n}$ if and only if $x_0 \in S$ and 
$${\rm rank} \left( \begin{array}{cc}
N_1 (x_0) & 0 \cr
\vdots & \vdots \cr
N_{c_s} (x_0) & 0 \cr
x_0 & t_0  \cr
 0 & 1 \cr
\end{array} \right) < c_S+2,$$
where $(N_1(x_0),\ldots,N_{c_S}(x_0))$ is a basis of the normal space to $S$ at $x_0$, that is if and only if 
$${\rm rank} \left( \begin{array}{c}
N_1 (x_0) \cr
\vdots  \cr
N_{c_s} (x_0)  \cr
x_0  \cr
\end{array} \right) < c_S+1.$$
But if $x_0$ is a point of $S$ close to $0$ but distinct from $0$, the sphere $\mathbb{S}_{\vert x_0 \vert}^{n-1}$ intersects the stratum $S$ transversally. We conclude that the unique possible critical point of $t_{\vert \widehat{S}\cap \mathbb{S}_{r_\epsilon}^n}$ is the point $(0,r_\epsilon)$. But $v$ is in $(-C_0 Y \setminus -C_0 Y')\cap \mathbb{S}^{n-1}$, so by Lemma \ref{Lemma4.14}, there is $\alpha >0$ such that $]0,\alpha] \cdot (-v) \subset T$ where $T$ is a stratum of $Y$ of dimension $n$. Hence $\{0\} \times ]0,\alpha]$ is included in $\widehat{T_v}$. We conclude that for $0< \epsilon \ll 1$, $(0,r_\epsilon)$ is the only critical point of $t_{\vert Z_v \cap \{ t > \frac{2\epsilon}{a}\} \cap \mathbb{B}_{r_\epsilon}^{n+1}}$. Moreover it is a strict local maximum. Applying Theorem 3.1 in \cite{DutertreManuscripta12} and Lemma \ref{Lemma4.15}, we get 
$$\chi \left( Z_v \cap \{t \ge \frac{2\epsilon}{a} \} \cap \mathbb{B}_{r_\epsilon}^{n+1} \right)
-\chi \left( Z_v \cap \{t = \frac{2\epsilon}{a} \} \cap \mathbb{B}_{r_\epsilon}^{n+1} \right)=0,$$
and 
$$\chi \left( Z_v \cap \{t \ge \frac{2\epsilon}{a} \} \cap \mathbb{B}_{r_\epsilon}^{n+1} \right) ={\rm ind}(-t, Z_v \cap \mathbb{B}_{r_\epsilon}^{n+1}, (0,r_\epsilon) )=1,$$
because $-t_{\vert Z_v \cap \mathbb{B}_{r_\epsilon}^{n+1}} $ has a strict local minimum at the point $(0,r_\epsilon)$. \endproof

\begin{corollary}\label{Corollary4.17}
Under the same assumptions, we have 
$$ \chi \left( {\rm Lk} (Z_v \cap \{ t \ge 0 \}) \right)=
\chi_c \left( Z_v \cap \{ t \ge 0 \} \cap \{ (x,t) \ \vert \
\vert x \vert=\epsilon, t < \frac{2\epsilon}{a} \}\right) + 1.$$
\end{corollary}
\proof Let $h(x,t)$ be the semi-algebraic function defined by 
$$h(x,t)={\rm max} \left( \vert x \vert, \frac{a}{2}t \right).$$
As explained by Durfee in \cite{Durfee}, Section 3, the link of $Z_v \cap \{ t \ge 0 \}$ is homeomorphic to $Z_v \cap \{ t \ge 0 \} \cap \{h=\epsilon\}$ for $0< \epsilon \ll 1$. We have 
$$Z_v \cap \{ t \ge 0 \} \cap \{h=\epsilon\} =
Z_v \cap \{ t \ge 0 \} \cap \{\vert x \vert= \epsilon, \frac{at}{2} \le \epsilon\} $$ $$\bigcup
Z_v \cap \{ t \ge 0 \} \cap \{\vert x \vert \le \epsilon, \frac{at}{2} = \epsilon\}.$$
It is enough to use the previous lemma and the additivity of $\chi_c$. \endproof

\subsection{On the two sides of the kinematic formulas}
We prove the existence of the left-hand sides of the kinematic formulas, and we show that both sides of the formulas are symmetric in $X$ and $Y$.
We also give a relation with the polar invariants. 

Let $(X,0) \subset (\mathbb{R}^n,0)$ and $(Y,0) \subset (\mathbb{R}^n,0)$ be two germs of closed definable sets. We assume that $X$ and $Y$ are included in an open set $U$. Let $\epsilon_0 >0$ be such that $\mathbb{B}_{\epsilon_0}^n \subset U$. 

1) Let us fix $(\epsilon,\delta)$ such that $0 \le \epsilon \le \epsilon_0$ and $0 \le \delta \le \epsilon$. Let 
$$\mathcal{A}= \left\{ (x,\gamma,v) \in \mathbb{R}^n \times SO(n) \times \mathbb{S}^{n-1} \ \vert \ 
x\in X, x-\delta v \in \gamma Y, \vert x \vert \le \epsilon \right\}.$$
It is a closed definable set. By Hardt's theorem applied to the projection $\pi : \mathcal{A} \to SO(n) \times \mathbb{S}^{n-1}$, the function $(\gamma, v) \mapsto  \chi \left( X \cap (\gamma Y +\delta v) \cap \mathbb{B}_\epsilon^n \right)$ takes a finite number of values. As in \cite{BernigBroeckerFourier} we equip $SO(n)$ with the Haar measure $d\gamma$, normalized in such a way that the volume of $SO(n)$ is $s_{n-1}$. We equip $\mathbb{S}^{n-1}$ with the usual Riemanniann measure (or density) $dv$ and $SO(n) \times \mathbb{S}^{n-1}$ with the product measure $d\gamma dv$. With this measure, the function $(\gamma,v) \mapsto \chi \left( X \cap (\gamma Y +\delta v) \cap \mathbb{B}_\epsilon^n \right)$ is integrable and so the integral
$$\int_{SO(n) \times \mathbb{S}^{n-1} } \chi \left( X \cap (\gamma Y +\delta v) \cap \mathbb{B}_\epsilon^n \right) d\gamma dv$$
exists and is finite. Moreover the function $$\gamma \mapsto 
\int_{\mathbb{S}^{n-1}} \chi \left( X \cap (\gamma Y +\delta v) \cap \mathbb{B}_\epsilon^n \right) dv $$ is integrable and 
the function $$ v \mapsto \int_{SO(n)} \chi \left( X \cap (\gamma Y +\delta v) \cap \mathbb{B}_\epsilon^n \right) d\gamma$$ is integrable and we have 
$$\displaylines{
\qquad  \int_{SO(n) \times \mathbb{S}^{n-1} } \chi \left( X \cap (\gamma Y +\delta v) \cap \mathbb{B}_\epsilon^n \right) d\gamma dv  \hfill \cr
\qquad \qquad \qquad=
\int_{SO(n)} [ \int_{ \mathbb{S}^{n-1} } \chi \left( X \cap (\gamma Y +\delta v) \cap \mathbb{B}_\epsilon^n \right) dv ] d\gamma \hfill \cr 
\hfill= \int_{\mathbb{S}^{n-1} } [ \int_{SO(n)} \chi \left( X \cap (\gamma Y +\delta v) \cap \mathbb{B}_\epsilon^n \right) d\gamma ] dv. \qquad \qquad \cr
}$$
Now let us fix $\epsilon >0$ such that $0 < \epsilon \le \epsilon_0$. By Hardt's theorem, for every $(\gamma,v) \in SO(n) \times \mathbb{S}^{n-1}$, there is a small interval $]0,\delta_\epsilon[$ such that the function $\delta \mapsto 
\chi \left( X \cap (\gamma Y +\delta v) \cap \mathbb{B}_\epsilon^n \right)$ is constant on $]0,\delta_\epsilon[$ and so $\lim_{\delta \to 0^+} \chi \left( X \cap (\gamma Y +\delta v) \cap \mathbb{B}_\epsilon^n \right)$ exists and is this constant value. Similarly as above, the function $(\gamma,v,\delta) \mapsto \chi \left( X \cap (\gamma Y +\delta v) \cap \mathbb{B}_\epsilon^n \right)$ takes a finite number of values and so, by Lebesgue's theorem, the function $(\gamma,v) \mapsto 
\lim_{\delta \to 0^+} \chi \left( X \cap (\gamma Y +\delta v) \cap \mathbb{B}_\epsilon^n \right)$ is integrable and we have
$$\displaylines{
\qquad \lim_{\delta \to 0^+} \int_{SO(n) \times \mathbb{S}^{n-1}} \chi \left( X \cap (\gamma Y +\delta v) \cap \mathbb{B}_\epsilon^n \right) d\gamma dv \hfill \cr
\hfill =
\int_{SO(n) \times \mathbb{S}^{n-1}} \lim_{\delta \to 0^+} \chi \left( X \cap (\gamma Y +\delta v) \cap \mathbb{B}_\epsilon^n \right) d\gamma dv . \qquad (*)\cr
}$$
Let us fix $(\gamma,v) \in SO(n) \times \mathbb{S}^{n-1}$ and let
$$\mathcal{B}= \Big\{ (x,\epsilon,\delta) \in \mathbb{R}^n \times \left\{ (\epsilon,\delta) \ \vert \
0 < \epsilon \le \epsilon_0 , 0 < \delta \le \epsilon \right\} \ \vert \ x\in X, x-\delta v \in Y, \vert x \vert \le \epsilon \Big\}.$$
It is a closed definable set. Applying the argument of the proof of Lemma \ref{Lemma3.1}, we see that there exists $0 < \epsilon_1 \le \epsilon_0$ and a definable function $r : ]0,\epsilon_1] \to \mathbb{R}$ continuous, monotone and strictly positive such that the function $$(\epsilon',\delta') \mapsto \chi \left( X \cap (\gamma Y +\delta v) \cap \mathbb{B}_{\epsilon'}^n \right)$$ is constant on $\{(\epsilon',\delta') \ \vert \ 0<\epsilon' < \epsilon_1, 0 < \delta' < r(\epsilon') \}$. But we see that for $\epsilon' \in ]0,\epsilon_1[$ and $0 < \delta' < r(\epsilon')$, 
$$\chi \left( X \cap (\gamma Y +\delta' v) \cap \mathbb{B}_{\epsilon'}^n \right) =
\lim_{\delta \to 0^+} \chi \left( X \cap (\gamma Y +\delta v) \cap \mathbb{B}_{\epsilon'}^n \right).$$
Therefore the limit $\lim_{\epsilon \to 0} \lim_{\delta \to 0^+} \chi \left( X \cap (\gamma Y +\delta v) \cap \mathbb{B}_\epsilon^n \right)$ exists and equals the above constant value. Always by Hardt's theorem, the function $(\gamma,v,\epsilon,\delta) \mapsto \chi \left( X \cap (\gamma Y +\delta v) \cap \mathbb{B}_{\epsilon}^n \right)$ takes a finite number of values and so does the function $(\gamma,v,\epsilon) \mapsto \lim_{\delta \to 0^+} \chi \left( X \cap (\gamma Y +\delta v) \cap \mathbb{B}_{\epsilon}^n \right)$. By Lebesgue's theorem, the function $(\gamma,v) \mapsto 
\lim_{\epsilon \to 0} \lim_{\delta \to 0^+} \chi \left( X \cap (\gamma Y +\delta v) \cap \mathbb{B}_{\epsilon}^n \right)$ is integrable and 
$$\lim_{\epsilon \to 0} \int_{SO(n) \times \mathbb{S}^{n-1}} 
\lim_{\delta \to 0^+} \chi \left( X \cap (\gamma Y +\delta v) \cap \mathbb{B}_{\epsilon}^n \right) d\gamma dv$$ $$ =
 \int_{SO(n) \times \mathbb{S}^{n-1}}  \lim_{\epsilon \to 0} \lim_{\delta \to 0^+}
 \chi \left( X \cap (\gamma Y +\delta v) \cap \mathbb{B}_{\epsilon}^n \right) d\gamma dv .$$
Finally, by Equality $(*)$ above, we have that $$\lim_{\epsilon \to 0} [\lim_{\delta \to 0^+}  \int_{SO(n) \times \mathbb{S}^{n-1}} 
\chi \left( X \cap (\gamma Y +\delta v) \cap \mathbb{B}_\epsilon^n \right) d\gamma dv]$$ exists and equals 
$$\int_{SO(n) \times \mathbb{S}^{n-1}}  \lim_{\epsilon \to 0} \lim_{\delta \to 0^+}
 \chi \left( X \cap (\gamma Y +\delta v) \cap \mathbb{B}_\epsilon^n \right) d\gamma dv.$$ 
\begin{definition}\label{Definition4.18}
{\rm For two germs of closed definable sets $(X,0)\subset (\mathbb{R}^n,0)$ and $(Y,0)\subset (\mathbb{R}^n,0)$, we set 
$$\sigma(X,Y,0) = \frac{1}{s_{n-1}^2} 
\int_{SO(n) \times \mathbb{S}^{n-1}}  \lim_{\epsilon \to 0} \lim_{\delta \to 0^+}
 \chi \left( X \cap (\gamma Y +\delta v) \cap \mathbb{B}_\epsilon^n \right) d\gamma dv.$$}
\end{definition}

2) Let us fix $\epsilon$ such that $0 \le \epsilon \le \epsilon_0$. Let
$$\mathcal{C}= \left\{ (x,\gamma,H) \in \mathbb{R}^n \times SO(n) \times G_n^{n-1} \ \vert \ x  \in X \cap \gamma Y \cap H, \vert x \vert = \epsilon \right\}.$$
It is a closed definable set. By Hardt's theorem applied to the projection $\pi : \mathcal{C} \to SO(n) \times G^n_{n-1}$, the function $(\gamma,H) \mapsto \chi \left( X \cap \gamma Y \cap H\cap \mathbb{S}_\epsilon^{n-1} \right)$ takes a finite number of values. As above, we deduce that the function $(\gamma,H) \mapsto \chi \left( {\rm Lk} (X \cap \gamma Y \cap H) \right)$ takes a finite number of values. We equip $G_n^{n-1}$ with the $SO(n)$-invariant measure (or density) $dH$  and $SO(n) \times G_n^{n-1}$ with the product measure. With this measure, the function $(\gamma,H) \mapsto \chi \left( {\rm Lk} (X \cap \gamma Y \cap H) \right)$ is integrable and so the integral $\int_{SO(n) \times G_n^{n-1}} \chi \left( {\rm Lk} (X \cap \gamma Y \cap H) \right) d\gamma dH$ exists. Moreover the function $\gamma \mapsto 
\int_{G_n^{n-1}} \chi \left( {\rm Lk} (X \cap \gamma Y \cap H) \right) dH$ is integrable and so is the function $$(\gamma,v) \mapsto \int_{G_n^{n-1}} \chi \left( {\rm Lk} (X \cap \gamma Y \cap H) \right) dH$$ on $SO(n) \times \mathbb{S}^{n-1}$. 
Similarly it is easy to see that the function $\gamma \mapsto \chi \left( {\rm Lk} (X \cap \gamma Y ) \right)$ is integrable on $SO(n)$ and so is the function $(\gamma,v) \mapsto \chi \left( {\rm Lk} (X \cap \gamma Y ) \right)$ on $SO(n) \times \mathbb{S}^{n-1}$.

3) By Lemmas \ref{Lemma4.10} and \ref{Lemma4.11}, there exists a definable subset $\Delta \subset SO(n) \times \mathbb{S}^{n-1}$ of positive codimension such that for $(\gamma,v) \notin \Delta$, 
$$\displaylines{
\quad \lim_{\epsilon \to 0} \lim_{\delta \to 0^+} 
\Lambda_0 \left(X \cap (\gamma Y + \delta v) ,  X \cap (\gamma Y + \delta v) \cap \mathbb{B}_\epsilon^n \right)= \hfill \cr 
\quad \quad \lim_{\epsilon \to 0} \lim_{\delta \to 0^+} \chi \left(X \cap (\gamma Y + \delta v) \cap \mathbb{B}_\epsilon^{n} \right)\hfill \cr -\frac{1}{2} \chi \left( {\rm Lk} (X\cap \gamma Y) \right)
-\frac{1}{2s_{n-1} }\int_{\mathbb{S}^{n-1}} \chi \left({\rm Lk} (X \cap \gamma Y \cap \{u^*=0\} ) \right)du.
}$$
Therefore the function $$(\gamma,v) \mapsto \lim_{\epsilon \to 0} \lim_{\delta \to 0^+} 
\Lambda_0 \left(X \cap (\gamma Y + \delta v) ,  X \cap (\gamma Y + \delta v) \cap \mathbb{B}_\epsilon^n \right)$$ is integrable on $SO(n) \times \mathbb{S}^{n-1}$. 
\begin{definition}\label{Definition4.19}
{\rm For two germs of closed definable set $(X,0)\subset (\mathbb{R}^n,0)$ and $(Y,0)\subset (\mathbb{R}^n,0)$, we set 
$$\displaylines{
\quad \Lambda_0^{{\rm lim}}(X,Y,0) \hfill \cr
\hfill = \frac{1}{s_{n-1}^2} 
\int_{SO(n) \times \mathbb{S}^{n-1}}  \lim_{\epsilon \to 0} \lim_{\delta \to 0^+}
\Lambda_0 \left(X \cap (\gamma Y + \delta v) ,  X \cap (\gamma Y + \delta v) \cap \mathbb{B}_\epsilon^n \right) d\gamma dv. \cr
}$$}
\end{definition}
We note that 
$$\displaylines{
\qquad \Lambda_0^{{\rm lim}}(X,Y,0) = \sigma(X,Y,0) -\frac{1}{2 s_{n-1}} \int_{SO(n)}\chi \left( {\rm Lk} (X\cap \gamma Y) \right) d\gamma \hfill \cr 
\hfill-\frac{1}{s_{n-1}^2 } \int_{SO(n)} \int_{\mathbb{S}^{n-1}} \chi \left({\rm Lk} (X \cap \gamma Y \cap \{u^*=0\}) \right)du d\gamma. \qquad \cr
}$$
The two limits $\sigma(X,Y,0)$ and $\Lambda_0^{{\rm lim}}(X,Y,0)$ are symmetric in $X$ and $Y$, as explained in the next proposition.
\begin{proposition}\label{Proposition4.20}
For two germs of closed definable sets $(X,0)\subset (\mathbb{R}^n,0)$ and $(Y,0)\subset (\mathbb{R}^n,0)$, we have 
$$\sigma(X,Y,0)=\sigma(Y,X,0) \hbox{ and } \Lambda_0^{{\rm lim}}(X,Y,0)=
\Lambda_0^{{\rm lim}}(Y,X,0).$$
\end{proposition}
\proof By Lemma \ref{Lemma4.6}, we know that 
$$\lim_{\epsilon \to 0} \lim_{\delta \to 0^+}
\chi \left( X \cap (\gamma Y + \delta v) \cap \mathbb{B}_\epsilon^n \right) = 
\lim_{\epsilon \to 0} \lim_{\delta \to 0^+}
\chi \left(( X - \delta v) \cap \gamma Y  \cap \mathbb{B}_\epsilon^n \right).$$
The change of variables $v \mapsto -v $ gives that 
$$\int_{\mathbb{S}^{n-1}} \chi \left( (X-\delta v) \cap \gamma Y  \cap \mathbb{B}_\epsilon^n \right) dv
 = \int_{\mathbb{S}^{n-1}} \chi \left( (X+\delta v) \cap \gamma Y  \cap \mathbb{B}_\epsilon^n \right) dv,$$
and so
$$\int_{\mathbb{S}^{n-1}} \chi \left( (X-\delta v) \cap \gamma Y  \cap \mathbb{B}_\epsilon^n \right) dv
 = \int_{\mathbb{S}^{n-1}} \chi \left( \gamma^{-1}(X+\delta v) \cap  Y  \cap \mathbb{B}_\epsilon^n \right) dv.$$
Hence, by the change of variables $\gamma \mapsto \gamma^{-1}$ on $SO(n)$, we are lead to compute 
$$\int_{SO(n)} [ \int_{\mathbb{S}^{n-1}} \chi \left( \gamma(X+\delta v) \cap  Y  \cap \mathbb{B}_\epsilon^n \right) dv ]d\gamma.$$
But for $\gamma \in SO(n)$, the change of variables $u \mapsto \gamma u$ gives that
$$\int_{\mathbb{S}^{n-1}} \chi \left( \gamma (X +\delta u) \cap Y \cap \mathbb{B}_\epsilon^n \right) du =
\int_{\mathbb{S}^{n-1}} \chi \left(( \gamma X +\delta v ) \cap Y \cap \mathbb{B}_\epsilon^n \right) dv.$$
Finally we get that 
$$\displaylines{
\qquad \int_{SO(n)} \int_{\mathbb{S}^{n-1}} \chi \left( \gamma (X +\delta v) \cap Y \cap \mathbb{B}_\epsilon^n \right)dv d\gamma = \hfill \cr
\hfill \int_{SO(n)} \int_{\mathbb{S}^{n-1}} \chi \left( (
\gamma X +\delta v) \cap Y \cap \mathbb{B}_\epsilon^n \right)dv d\gamma .
\qquad \cr
}$$
It is enough to pass to the limits to get the equality $\sigma(X,Y,0)=\sigma(Y,X,0)$. The result for $\Lambda_0^{\rm lim}(X,Y,0)$ is obtained applying the relation between 
$\sigma(X,Y,0)$ and $\Lambda_0^{\rm lim}(X,Y,0)$. \endproof

Now let us relate $\sigma(X,Y,0)$ with the polar invariants of Comte and Merle \cite{ComteMerle}. 

\begin{proposition}\label{Proposition4.21}
If $H \in G_n^{n-k}$, $k \in \{0,\ldots,n\}$, then we have
$$\sigma(X,H,0)= \sigma_k(X,0).$$
\end{proposition}
\proof The proof is straightforward for $k=0$, because in this case $\sigma(X,\mathbb{R}^n,0)=\lim_{\epsilon \to 0} \chi (X \cap \mathbb{B}_\epsilon^n)=1$. 

Let us assume that $k >0$. First we note that
$$\sigma(X,H,0)= \frac{1}{g_n^{n-k}} \int_{G_n^{n-k}} \left( \frac{1}{s_{n-1}}
\int_{\mathbb{S}^{n-1}} \lim_{\epsilon \to 0} \lim_{\delta \to 0^+} 
\chi \left( X \cap (H +\delta v) \cap \mathbb{B}_\epsilon^n \right) dv \right) dH.$$ 
For $H \in G_n^{n-k}$, we recall that $\mathbb{S}_{H^\perp}^{k-1}$ is the unit sphere in $H^\perp$ and we denote by $p_{H^\perp}$ the orthogonal projection onto $H^\perp$. If $v \in \mathbb{S}_{H^\perp}^{k-1}$ and $w \in \mathbb{S}^{n-1}$ are such that $\frac{p_{H^\perp}(w)}{\vert p_{H^\perp}(w) \vert }=v$, then
$$\lim_{\epsilon \to 0} \lim_{\delta \to 0^+} \chi \left( X \cap (H+\delta v) \cap \mathbb{B}_\epsilon^n \right)= \lim_{\epsilon \to 0} \lim_{\delta \to 0^+} \chi \left( X \cap (H+\delta w) \cap \mathbb{B}_\epsilon^n \right). $$
This implies that 
$$\sigma(X,H,0) = \frac{1}{g_n^{n-k}} \int_{G_n^{n-k}} \left( \frac{1}{s_{k-1}} \int_{\mathbb{S}_{H^\perp}^{k-1}} \lim_{\epsilon \to 0} \lim_{\delta \to 0^+} \chi \left( X \cap (H+\delta v) \cap \mathbb{B}_\epsilon^n \right) dv \right) dH,$$
(see the proof of Corollary 5.7 in \cite{DutertreAdvGeo2008} for a similar argument or use the co-area formula). To end the proof, it is enough to show, with the notations of Section 3, that
$$ \sum_{i=1}^{N_P} \chi_i^P \cdot \Theta(K_i^p,0) = 
\frac{1}{s_{k-1}} \int_{\mathbb{S}_P^{k-1}} 
\lim_{\epsilon \to 0} \lim_{\delta \to 0^+} \chi \left( X \cap (P^\perp+\delta v) \cap \mathbb{B}_\epsilon^n \right) dv .$$
By Lemma \ref{Lemma4.14}, we know that if $v \in C_0 K_i^P \setminus C_0( \overline{K_i^P} \setminus K_i^P )$, there is $\delta >0$ such that $]0,\delta]\cdot v \subset K_i^P$. Hence 
$$\frac{1}{s_{k-1}} \int_{\mathbb{S}_P^{k-1}} \lim_{\epsilon \to 0} \lim_{\delta \to 0^+} \chi \left( X \cap (P^\perp+\delta v) \cap \mathbb{B}_\epsilon^n \right) dv   =
\sum_{i=1}^{N_P} \chi_i^P \cdot \frac{{\rm vol} ( C_0 K_i^P \cap \mathbb{S}_P^{k-1})}{s_{k-1}}.$$
By \cite{KurdykaRaby} Lemma 2.1, $\frac{{\rm vol} ( C_0 K_i^P \cap \mathbb{S}_P^{k-1})}{s_{k-1}}$ is exactly $\Theta_k (K_i^P,0)$. \endproof

\begin{remark}
{\rm The equality $$\sigma_k(X,0)= \frac{1}{g_n^{n-k}} \int_{G_n^{n-k}} \left( \frac{1}{s_{k-1}} \int_{\mathbb{S}_{H^\perp}^{k-1}} \lim_{\epsilon \to 0} \lim_{\delta \to 0^+} \chi \left( X \cap (H+\delta v) \cap \mathbb{B}_\epsilon^n \right) dv \right) dH,$$
is natural because the two sides of the equality coincide for they measure the same mean-value of Euler characteristics. It
already appeared in \cite{ComteMerle} page 244 in the conic case, and we used it in \cite{DutertreJofSingProcTrot,DutertreIsrael,DutertreAdvGeo2019} as a definition for the polar invariants. We prove it here for completeness.}
\end{remark}
We end this subsection  with another symmetry result.
\begin{lemma}\label{Lemma4.23}
For two germs of closed definable sets $(X,0)\subset (\mathbb{R}^n,0)$ and $(Y,0)\subset (\mathbb{R}^n,0)$, we have 
$$\sum_{i=0}^n \Lambda_i^{\rm lim} (X,0)\cdot \sigma_{n-i}(Y,0) =
\sum_{i=0}^n \sigma_i (X,0)\cdot \Lambda_{n-i}^{\rm lim}(Y,0).$$
\end{lemma}
\proof By Theorem \ref{Theorem3.8}, we get 
$$\displaylines{\qquad  \sum_{i=0}^n \Lambda_i^{\rm lim} (X,0)\cdot \sigma_{n-i}(Y,0) = 
\sum_{i=0}^{n-1} \left(  \sigma_i (X,0) - \sigma_{i+1}(X,0) \right) \cdot \sigma_{n-i}(Y,0) \hfill \cr
\qquad  \qquad  + \sigma_n(X,0) \cdot \sigma_0 (Y,0) 
=\sum_{i=0}^{n-1} \sigma_i (X,0)  \cdot \sigma_{n-i}(Y,0) \hfill \cr 
\hfill -\sum_{i=0}^{n-1} \sigma_{i+1} (X,0)  \cdot \sigma_{n-i}(Y,0)
+ \sigma_n(X,0) \cdot \sigma_0 (Y,0). \qquad \qquad 
}$$
Therefore we obtain
$$\displaylines{
\quad \sum_{i=0}^n \Lambda_i^{\rm lim} (X,0)\cdot \sigma_{n-i}(Y,0) 
= \sigma_0(X,0) \cdot \sigma_n (Y,0) \hfill \cr
\hfill  +
\sum_{i=1}^{n-1}   \sigma_i (X,0)\cdot \left(  \sigma_{n-i}(Y,0)-\sigma_{n-i+1}(Y,0) \right) +
\sigma_n (X,0) \cdot \left(  \sigma_0(Y,0)-\sigma_1(Y,0) \right) \quad \cr
}$$
Another application of Theorem \ref{Theorem3.8} gives the result.\endproof

\section{A new spherical kinematic formula}
We give a new spherical kinematic formula for two definable subsets of the unit sphere. 

Let $X \subset \mathbb{S}^{n-1}$ be a compact definable set and let $Y \subset \mathbb{S}^{n-1}$ be a definable set, not necessarly compact. We recall that the $\tilde{\Lambda}_i$'s, $i=0,\ldots,n-1$, denote the spherical Lipschitz-Killing measures.
\begin{proposition}\label{Proposition5.1}
The following kinematic formula holds:
$$\frac{1}{s_{n-1}}\int_{SO(n)} \chi_c (X \cap \gamma Y) d\gamma =
\sum_{i=0}^{n-1} \frac{\tilde{\Lambda}_i(X,X)}{s_i} \cdot \frac{1}{g_n^{i+1}} \int_{G_n^{i+1}} \chi_c (Y \cap H) dH.$$
\end{proposition}
\proof {\bf First step:} We study the case $Y$ compact. Applying the generalized spherical  Gauss-Bonnet formula (see Theorem 1.2 in \cite{BernigBroeckerFourier}) to $\chi(X \cap \gamma Y)$, we obtain
$$\int_{SO(n)} \chi(X \cap \gamma Y) d\gamma = \sum_{i=0,2,\ldots}^{n-1}
\frac{2}{s_i} \int_{SO(n)} \tilde{\Lambda}_i (X \cap \gamma Y,X \cap \gamma Y) d\gamma.$$
Then we apply the generalized spherical kinematic formula (see \cite{BernigBroeckerFourier,FuIndiana90}) to each $\tilde{\Lambda}_i(X \cap \gamma Y, X \cap \gamma Y)$ and we get 
$$\int_{SO(n)} \chi(X \cap \gamma Y) d \gamma =
\sum_{i=0,2,\ldots}^{n-1} \frac{2}{s_i} \sum_{p+q=i+n-1} \frac{s_i s_{n-1}}{s_p s_q}
\tilde{\Lambda}_p(X,X) \tilde{\Lambda}_q (Y,Y).$$
Therefore we have 
$$\displaylines{
\quad \frac{1}{s_{n-1}} \int_{SO(n)} \chi( X \cap \gamma Y) d \gamma = 
\sum_{i=0,2,\ldots}^{n-1} \sum_{p+q=i+n-1} \frac{\tilde{\Lambda}_p (X,X)}{s_p}
\frac{2\tilde{\Lambda}_q (Y,Y)}{s_q}\hfill \cr
\quad \quad = \sum_{p+q=n-1} \frac{\tilde{\Lambda}_p (X,X)}{s_p}
\frac{2\tilde{\Lambda}_q (Y,Y)}{s_q}+\sum_{p+q=n+1} \frac{\tilde{\Lambda}_p (X,X)}{s_p}
\frac{2\tilde{\Lambda}_q (Y,Y)}{s_q}\hfill \cr 
\hfill +\cdots + \sum_{p+q=e(n)} \frac{\tilde{\Lambda}_p (X,X)}{s_p}\frac{2\tilde{\Lambda}_q (Y,Y)}{s_q},\quad
}$$
where $e(n)=2n-2$ if $n-1$ is even or $2n-3$ if $n-1$ is odd. This equality can be rewritten in the following way:
$$\displaylines{
\quad \frac{1}{s_{n-1}} \int_{SO(n)} \chi( X \cap \gamma Y) d \gamma =  
\frac{\tilde{\Lambda}_{n-1} (X,X)}{s_{n-1}}\left( \sum_{q=0,2,\ldots} \frac{2\tilde{\Lambda}_q (Y,Y)}{s_q} \right) \hfill \cr
\hfill + \frac{\tilde{\Lambda}_{n-2} (X,X)}{s_{n-2}}\left( \sum_{q=1,3,\ldots} \frac{2\tilde{\Lambda}_q (Y,Y)}{s_q} \right)
+\cdots+ \frac{\tilde{\Lambda}_0 (X,X)}{2}\left(  \frac{2\tilde{\Lambda}_{n-1} (Y,Y)}{s_{n-1}} \right). \quad
}$$
In \cite{DutertreGeoDedicata2012} pages 175-176, we proved that
$$\frac{1}{g_n^1} \int_{G_n^1} \chi (Y \cap H) dH = \frac{2\tilde{\Lambda}_{n-1} (Y,Y)}{s_{n-1}},\
\frac{1}{g_n^2} \int_{G_n^2} \chi (Y \cap H) dH = \frac{2\tilde{\Lambda}_{n-2} (Y,Y)}{s_{n-2}},$$
and for $k \ge 3$,
$$\displaylines{
\qquad \frac{1}{g_n^{k-2}} \int_{G_n^{k-2}} \chi (Y \cap H) dH =
\sum_{i=2,4,\ldots} \frac{2 \tilde{\Lambda}_{n-k+i} (Y,Y)}{s_{n-k+i}} \hfill \cr
\hfill = \sum_{q=n-k+2,n-k+4,\ldots} \frac{2 \tilde{\Lambda}_{q} (Y,Y)}{s_{q}}. \qquad \cr
}$$
Applying these relations, we get that
$$\frac{1}{s_{n-1}} \int_{SO(n)} \chi ( X \cap \gamma Y) d \gamma =
\frac{\tilde{\Lambda}_{n-1} (X,X)}{s_{n-1}} \chi(Y)+ \frac{\tilde{\Lambda}_{n-2} (X,X)}{s_{n-2}} \frac{1}{g_n^{n-1}}\int_{G_n^{n-1}} \chi (Y \cap H) dH$$
$$+\frac{\tilde{\Lambda}_{n-3} (X,X)}{s_{n-3}} \frac{1}{g_n^{n-2}}\int_{G_n^{n-2}} \chi (Y \cap H) dH + \cdots + \frac{\tilde{\Lambda}_0 (X,X)}{s_0} \frac{1}{g_n^1}\int_{G_n^1} \chi (Y \cap H) dH.$$

{\bf Second step:} Let $Y \subset \mathbb{S}^{n-1}$ be compact and let $K \subsetneq Y$ be a compact definable set. By the first step, we have
$$\frac{1}{s_{n-1}}\int_{SO(n)} \chi (X \cap \gamma Y) d\gamma =
\sum_{i=0}^{n-1} \frac{\tilde{\Lambda}(X,X)}{s_i} \cdot \frac{1}{g_n^{i+1}} \int_{G_n^{i+1}} \chi (Y \cap H) dH,$$
and
$$\frac{1}{s_{n-1}}\int_{SO(n)} \chi (X \cap \gamma K) d\gamma =
\sum_{i=0}^{n-1} \frac{\tilde{\Lambda}(X,X)}{s_i} \cdot \frac{1}{g_n^{i+1}} \int_{G_n^{i+1}} \chi (K \cap H) dH.$$
For each $\gamma \in SO(n)$, $\gamma Y = (\gamma Y \setminus \gamma K) \sqcup \gamma K = \gamma (Y \setminus K) \sqcup \gamma K$ because $\gamma$ is bijective. Hence $\chi (X \cap \gamma Y)=\chi_c \left( X \cap \gamma (Y \setminus K)\right) + \chi (X \cap \gamma K) $ and
$$\displaylines{
\qquad  \frac{1}{s_{n-1}} \int_{SO(n)} \chi_c \left( X \cap \gamma (Y \setminus K) \right) d \gamma \hfill \cr
\qquad \qquad =\sum_{i=0}^{n-1} \frac{\tilde{\Lambda}_i(X,X)}{s_i} \cdot
\frac{1}{g_n^{i+1}} \int_{G_n^{i+1}} \left[ \chi (Y \cap H) -\chi(K \cap H) \right] dH \hfill \cr
\hfill = \sum_{i=0}^{n-1} \frac{\tilde{\Lambda}_i(X,X)}{s_i} \cdot
\frac{1}{g_n^{i+1}} \int_{G_n^{i+1}}  \chi_c \left( (Y\setminus K) \cap H \right) dH. \qquad \cr
}$$
This gives the result for the set $Y \setminus K$.

{\bf Third step:} We prove the general case. Since $Y$ is definable, it admits the following cell decomposition $Y=\sqcup_{j=1}^r C_j$, where each $C_j$ is a definable subset homeomorphic to a unit cube $]0,1[^{d_j}$. By the second step, the result is valid for each cell $C_j$, because $\overline{C_j}$ and $\overline{C_j} \setminus C_j$ are compact and definable. By additivity of $\chi_c$, we have 
$$\displaylines{
\quad \frac{1}{s_{n-1}} \int_{SO(n)} \chi_c (X \cap \gamma Y) d\gamma =
\frac{1}{s_{n-1}} \int_{SO(n)} \chi_c \left(X \cap \gamma (\sqcup_{j=1}^r C_j) \right) d\gamma \hfill \cr 
\quad \quad = \frac{1}{s_{n-1}} \int_{SO(n)} \chi_c \left(X \cap  (\sqcup_{j=1}^r \gamma C_j) \right) d\gamma  = 
\frac{1}{s_{n-1}} \int_{SO(n)} \chi_c \left(\sqcup_{j=1}^r (X \cap \gamma  C_j )\right) d\gamma \hfill \cr
 = \sum_{j=1}^r \frac{1}{s_{n-1}} \int_{SO(n)} \chi_c (X \cap \gamma C_j) d\gamma . \cr
}$$
Applying the second step, we obtain 
$$\displaylines{
\quad \frac{1}{s_{n-1}} \int_{SO(n)} \chi_c (X \cap \gamma Y) d\gamma =
\sum_{j=1}^r \sum_{i=0}^{n-1} \frac{\tilde{\Lambda}_i(X,X)}{s_i} \cdot
\frac{1}{g_n^{i+1}} \int_{G_n^{i+1}}  \chi_c ( C_j \cap H) dH \hfill \cr 
\qquad \qquad = \sum_{i=0}^{n-1} \frac{\tilde{\Lambda}_i(X,X)}{s_i} \cdot
\frac{1}{g_n^{i+1}} \int_{G_n^{i+1}} \sum_{j=1}^r \chi_c (C_j \cap H) dH \hfill \cr 
\hfill =\sum_{i=0}^{n-1} \frac{\tilde{\Lambda}_i(X,X)}{s_i} \cdot
\frac{1}{g_n^{i+1}} \int_{G_n^{i+1}} \chi_c (Y \cap H) dH. \qquad \cr
}$$ \endproof

\section{A second kinematic formula in the unit ball}
We deduce from the previous spherical kinematic a new kinematic formula for definable subsets of the unit ball. 

Let $X \subset \mathbb{R}^n$ be a closed conic definable set. Let $Y \subset \mathbb{B}^n$ be another definable set.
\begin{proposition}\label{Proposition6.1} The following kinematic formula holds:
$$\frac{1}{s_{n-1}} \int_{SO(n)} \chi_c ( X \cap \gamma Y) d\gamma = \sum_{i=0}^n \frac{\Lambda_i (X,X \cap \mathbb{B}^n)}{b_i} \cdot \frac{1}{g_n^i}\int_{G_n^i} \chi_c (Y \cap H) dH.$$
\end{proposition}
\proof Let us assume first that $0 \notin Y$ and let $\phi$ be the following definable mapping:
$$\begin{array}{ccccc}
\phi & : & \mathbb{B}^n \setminus \{0\} & \to & \mathbb{S}^{n-1} \cr
 &  & x & \mapsto & \frac{x}{\vert x \vert}. \cr
\end{array}$$
By Hardt's theorem, there exists a definable partition of $\phi (Y)$, $\phi (Y)= \sqcup_{j=1}^r W_j$, such that for $j \in \{1,\ldots,r\}$, the mapping $\phi_{\vert Y \cap \phi^{-1} (W_j) } :  \phi^{-1}(W_j) \cap Y \to W_j$ is trivial. By additivity and multiplicity of $\chi_c$, we can write
$\chi_c (Y)= \sum_{j=1}^r \alpha_j \chi_c(W_j)$, where $\alpha_j = \chi_c (F_j)$ with $F_j$ the fibre of $\phi_{\vert Y \cap \phi^{-1} (W_j) }$. Let us set $X^* = (X \cap \mathbb{B}^n) \setminus \{0\}$. 

If $w$ belongs to $\phi (Y) \cap \phi(X^*)$ then $w=\phi (y)=\phi (x)$ with $y \in Y$ and $x \in X^*$. Since $X$ is conic, $y$ belongs to $X^*$ and so $\phi (Y \cap X^*)=\phi (Y) \cap \phi (X^*)$. Therefore $\phi (Y \cap X^*)= \sqcup_{j=1}^r W_j \cap \phi(X^*)$. Note that if $w \in \phi(X^*)$ then $\phi^{-1}(w) \subset X^*$ by the conic structure of $X$. Hence if $w_j \in W_j \cap \phi(X^*)$, $\phi^{-1}(w_j) \cap Y \cap X =\phi^{-1}(w_j) \cap Y$ and $\chi_c ( \phi^{-1}(w_j) \cap Y \cap X)= \alpha_j$. Applying again Hardt's theorem to each $W_j \cap \phi (X^*)$ if necessary, one can conclude as above that 
$$\chi_c (Y \cap X)=\sum_{j=1}^r \alpha_j \chi_c (W_j \cap \phi(X^*) ).$$
Let $\gamma \in SO(n)$. Since $\phi \circ \gamma = \gamma \circ \phi$ and $\gamma$ is a definable homeomorphism, $\phi (\gamma Y)= \sqcup_{j=1}^r \gamma W_j$ and for $j \in \{1,\ldots,r\}$, the mapping 
$$\phi_{\vert \gamma Y \cap \phi^{-1}(\gamma W_j)} :  \gamma Y  \cap \phi^{-1}\left(\gamma W_j \right) \to \gamma W_j$$ is trivial, with fibre homeomorphic to $F_j$. As above, we can write 
$$\chi_c \left(X \cap \gamma Y \right)=\sum_{j=1}^r \alpha_j \chi_c \left(\gamma W_j \cap \phi (X^*) \right).$$
We can apply Proposition \ref{Proposition5.1} to the sets $W_j$ and $\phi(X^*)$. We get 
$$\displaylines{
\quad \frac{1}{s_{n-1}} \int_{SO(n)} \chi_c \left( X \cap \gamma Y \right) d \gamma =
\sum_{j=1}^r \alpha_j \frac{1}{s_{n-1}} \int_{SO(n)} \chi_c (\gamma W_j \cap \phi(X^*))d\gamma \hfill \cr
\quad \quad  =\sum_{j=1}^r \alpha_j 
\left( \sum_{i=1}^{n-1} \frac{\tilde{\Lambda}_i(\phi(X^*),\phi(X^*) )}{s_i} \cdot
\frac{1}{g_n^{i+1}}\int_{G_n^{i+1}} \chi_c (W_j \cap H) dH \right) \hfill \cr
 = \sum_{i=1}^{n-1} \frac{\tilde{\Lambda}_i(\phi(X^*),\phi(X^*) )}{s_i} \cdot
\frac{1}{g_n^{i+1}}\int_{G_n^{i+1}}\sum_{j=1}^r \alpha_j \chi_c (W_j \cap H) dH .  \cr
}$$
Since $H$ is conic, $\sum_{j=1}^r \alpha_j \chi_c (W_j \cap H)= \chi_c (Y \cap H)$ by the above argument. Applying Corollary 3.5 in \cite{DutertreGeoDedicata2012}, we obtain the following equality:
$$\displaylines{
\quad \frac{1}{s_{n-1}} \int_{SO(n)} \chi_c ( X \cap \gamma Y) d\gamma = \sum_{i=0}^{n-1} \frac{\Lambda_{i+1} (X,X\cap \mathbb{B}^n)}{b_{i+1}} \cdot \frac{1}{g_n^{i+1}}\int_{G_n^{i+1}} \chi_c (Y \cap H) dH \hfill \cr
 =\sum_{i=1}^{n} \frac{\Lambda_{i} (X,X\cap \mathbb{B}^n)}{b_{i}}\cdot \frac{1}{g_n^i}\int_{G_n^{i}} \chi_c (Y \cap H) dH,  \cr
}$$
which is the expected one when $0 \notin Y$. 

If $0 \in Y$ then $\chi_c(X \cap \gamma Y)= \chi_c(X \cap \gamma Y^*)+1$ and $\chi_c (Y \cap H)= \chi_c (Y^* \cap H)+1$ where $Y^*=Y \setminus \{0\}$. Therefore we have 
$$\displaylines{
\qquad \frac{1}{s_{n-1}} \int_{SO(n)} \chi_c ( X \cap \gamma Y) d\gamma =
1 + \frac{1}{s_{n-1}} \int_{SO(n)} \chi_c ( X \cap \gamma Y^*) d\gamma \hfill \cr
\qquad \qquad = 1 + \sum_{i=1}^{n} \frac{\Lambda_{i} (X,X\cap \mathbb{B}^n)}{b_{i}} \cdot\frac{1}{g_n^i}\int_{G_n^{i}} \chi_c (Y^* \cap H) dH \hfill \cr
\hfill = 1 + \sum_{i=1}^{n} \frac{\Lambda_{i} (X,X\cap \mathbb{B}^n)}{b_{i}} \cdot \frac{1}{g_n^i}\int_{G_n^{i}} \chi_c (Y \cap H) dH - \sum_{i=1}^{n} \frac{\Lambda_{i} (X,X \cap \mathbb{B}^n)}{b_{i}}. \qquad \cr
}$$
But by Corollary 5.2 in \cite{DutertreGeoDedicata2012}, we know that $1-\sum_{i=1}^{n} \frac{\Lambda_{i} (X,X\cap \mathbb{B}^n)}{b_{i}} = \Lambda_0(X,X\cap \mathbb{B}^n).$  \endproof

\section{A infinitesimal kinematic formula for conic sets}
We prove a first version of the infinitesimal principal kinematic formula for closed conic definable sets.

Let $X,Y \subset \mathbb{R}^n$ be two closed conic definable sets. We keep the notations used in Section 4.
\begin{lemma}\label{Lemma7.1}
There exists $\delta_1 >0$ such that for $0 < \delta \le \delta_1$, the topological types of $Z_{v,\delta} \cap \mathbb{B}^n$ does not depend on the choice $\delta$. Moreover, we have 
$$\lim_{\delta \to 0} \chi (Z_{v,\delta} \cap \mathbb{B}^n) = \lim_{\epsilon \to 0}
 \lim_{\delta \to 0} \chi ( Z_{v,\delta} \cap \mathbb{B}_\epsilon^n ).$$
\end{lemma}
\proof By Lemma \ref{Lemma4.4}, we know that there exists $\epsilon_0 >0$ such that for $0< \epsilon \le \epsilon_0$, there exits $\delta_\epsilon$ such that for $0< \delta \le \delta_\epsilon$, the topological type of $Z_{v,\delta} \cap \mathbb{B}_\epsilon^n$ does not depend on the choice of the couple $(\epsilon,\delta)$. Let us fix such a couple $(\epsilon,\delta)$. Let $\theta_\epsilon : \mathbb{R}^n \to \mathbb{R}^n$ be the diffeomorphism $\theta_\epsilon(x)=\frac{1}{\epsilon}x$. Then $\theta_\epsilon(Z_{v,\delta} \cap \mathbb{B}_\epsilon^n)= Z_{v,\frac{\delta}{\epsilon}} \cap \mathbb{B}^n$. Since $\lim_{\delta \to 0} \frac{\delta}{\epsilon}=0$, we get the result. \endproof

We are in position to state a first infinitesimal kinematic formula in the conic setting.
\begin{proposition}\label{Proposition7.2}
Let $X,Y \subset \mathbb{R}^n$ be two closed conic definable sets. The following kinematic formula holds:
$$\displaylines{
\qquad \frac{1}{s_{n-1}^2} \int_{SO(n) \times \mathbb{S}^{n-1}} \lim_{\delta \to 0^+}
\chi \left( X \cap (\gamma Y +\delta v) \cap \mathbb{B}^n \right)  d\gamma dv \hfill \cr \hfill =
\sum_{i=0}^n \frac{\Lambda_i ( X ,X \cap \mathbb{B}^n)}{b_i} \cdot
\sigma_{n-i}(Y,0). \qquad \cr
}$$
\end{proposition}
\proof Let us fix $\delta >0$. By the change of variable $u=\gamma v$, we have that for $\gamma \in SO(n)$
$$\int_{\mathbb{S}^{n-1}} \chi \left( X \cap (\gamma Y +\delta u) \cap \mathbb{B}^n \right) du 
=\int_{\mathbb{S}^{n-1}} \chi \left( X \cap \gamma (Y +\delta v) \cap \mathbb{B}^n \right) dv .$$
Applying Proposition \ref{Proposition6.1} to $X \cap \mathbb{B}^n$ and $(Y+\delta v) \cap \mathbb{B}^n$, we get that 
$$\displaylines{
\quad \frac{1}{s_{n-1}^2} \int_{SO(n) \times \mathbb{S}^{n-1}}
\chi \left( X \cap \gamma (Y +\delta v) \cap \mathbb{B}^n \right) d\gamma dv \hfill \cr 
\qquad \qquad = \frac{1}{s_{n-1}^2}  \int_{\mathbb{S}^{n-1}} \int_{SO(n)} 
\chi \left( X \cap \gamma (Y +\delta v) \cap \mathbb{B}^n \right) d\gamma dv \hfill \cr 
=\sum_{i=1}^n \frac{\Lambda_i ( X ,X \cap \mathbb{B}^n)}{b_i} \cdot \frac{1}{g_n^i s_{n-1}}
\int_{\mathbb{S}^{n-1}} \int_{G_n^i}   \chi \left( (Y+\delta v) \cap \mathbb{B}^n \cap H \right) dH dv ,\cr
}$$
and so that 
$$\displaylines{
\quad \frac{1}{s_{n-1}^2} \int_{SO(n) \times \mathbb{S}^{n-1}}
\chi \left( X \cap \gamma (Y +\delta v) \cap \mathbb{B}^n \right) dv d\gamma \hfill \cr  =\sum_{i=1}^n \frac{\Lambda_i ( X ,X \cap \mathbb{B}^n)}{b_i} \cdot \frac{1}{g_n^i}
\int_{G_n^i}  \frac{1}{s_{n-1}} \int_{\mathbb{S}^{n-1}} \chi \left( (Y+\delta v) \cap \mathbb{B}^n \cap H \right) dv  dH.
}$$
Passing to the limit as $\delta \to 0^+$ and using Lebesgue's theorem, we obtain that 
$$\displaylines{
\quad \frac{1}{s_{n-1}^2} \int_{SO(n) \times \mathbb{S}^{n-1}} 
\lim_{\delta \to 0^+} \chi \left( X \cap (\gamma Y +\delta v) \cap \mathbb{B}^n \right) dv d\gamma \hfill \cr
= \sum_{i=1}^n \frac{\Lambda_i ( X ,X \cap \mathbb{B}^n)}{b_i} \cdot \frac{1}{g_n^i}
\int_{G_n^i}  \frac{1}{s_{n-1}} \int_{\mathbb{S}^{n-1}} \lim_{\delta \to 0^+} \chi \left( (Y+\delta v) \cap \mathbb{B}^n \cap H \right) dv  dH.
}$$
By Lemmas \ref{Lemma4.6} and \ref{Lemma7.1} applied to $Y$ and $H$, we have that
$$\lim_{\delta \to 0^+} \chi \left( (Y+\delta v) \cap \mathbb{B}^n \cap H \right)=
 \lim_{\delta \to 0^+} \chi \left( Y \cap \mathbb{B}^n \cap (H-\delta v) \right),$$
and so 
$$\frac{1}{g_n^i} \int_{G_n^i} \frac{1}{s_{n-1}} \int_{\mathbb{S}^{n-1}}
\lim_{\delta \to 0^+} \chi \left( (Y+\delta v) \cap \mathbb{B}^n \cap H \right) dv dH =
\sigma_{n-i} (Y,0),$$
by Proposition \ref{Proposition4.21}. \endproof

\section{The principal kinematic formulas}
We prove our main results : the principal kinematic formulas for germs of closed definable sets. We will use the kinematic formula for closed conic definable sets proved in the previous section. We will proceed in several steps.

We keep the notations used in Section 4. For convenience, we also use the notation $\omega(x)$ for $\vert x \vert$, if $x$ is in $\mathbb{R}^n$.

{\bf First step:} $(X,0) \subset (\mathbb{R}^n,0)$ is a germ of closed definable set, $Y \subset \mathbb{R}^n$ is a closed conic definable set.

We assume that $X$ is included in an open neighborhood $U$ of $0$. Let $\epsilon_0 >0$ be such that $\mathbb{B}^n_{\epsilon_0} \subset U$. For $0< u \le \epsilon_0$, we set $X_u= X \cap \mathbb{S}_u^{n-1}$ and we denote by $CX_u$ the cone over $X_u$, i.e.:
$$CX_u = \left\{ x \in \mathbb{R}^n \ \vert \ \exists \lambda \in \mathbb{R}^+ \hbox{ and } z \in X_u 
\hbox{ such that } x = \lambda z \right\}.$$
\begin{lemma}\label{Lemma8.1}
There exists a definable subset $\Delta_Y \subset \mathbb{S}^{n-1}$ of positive codimension such that for $v \notin \Delta_Y$, 
$$\lim_{u \to 0} \chi \left( {\rm Lk} ( \widehat{CX_u} \cap \widehat{Y_v} \cap \{t \ge 0 \} )\right) =
\chi \left( {\rm Lk} ( \widehat{X} \cap \widehat{Y_v} \cap \{t \ge 0 \}) \right).$$
\end{lemma}
\proof If $v \notin (-Y) \cap \mathbb{S}^{n-1}$, then by Corollary \ref{Corollary4.13}, we have
$$\chi \left( {\rm Lk} ( \widehat{CX_u} \cap \widehat{Y_v} \cap \{t \ge 0 \} )\right) =
\lim_{\epsilon \to 0} \chi \left(\widehat{CX_u} \cap \widehat{Y_v} \cap \{t \ge 0 \} \cap \{ \omega=\epsilon\} \right)$$
and 
$$\chi \left( {\rm Lk} ( \widehat{X} \cap \widehat{Y_v} \cap \{t \ge 0 \} )\right) =
\lim_{\epsilon \to 0} \chi \left(\widehat{X} \cap \widehat{Y_v} \cap \{t \ge 0 \} \cap \{ \omega=\epsilon\} \right).$$
Let us choose $\epsilon_1 \ge 0$ such that for $0< u \le \epsilon_1$,
$$\chi \left( {\rm Lk} ( \widehat{X} \cap \widehat{Y_v} \cap \{t \ge 0 \} ) \right) =
\chi \left(\widehat{X} \cap \widehat{Y_v} \cap \{t \ge 0 \} \cap \{ \omega=u \} \right).$$
But $\widehat{X} \cap \{\omega = u \} = \widehat{CX_u}\cap \{ \omega = u\}$ and so 
$$\chi \left( \widehat{X} \cap \widehat{Y_v} \cap \{t \ge 0 \} \cap \{ \omega=u \} \right) =
\chi \left( \widehat{CX_u} \cap \widehat{Y_v} \cap \{t \ge 0 \} \cap \{ \omega=u \} \right).$$
For any $\epsilon > 0$, the mapping 
$$\begin{array}{ccc}
\widehat{CX_u} \cap \widehat{Y_v} \cap \{t \ge 0 \} \cap \{ \omega=\epsilon \} & \to & \widehat{CX_u} \cap \widehat{Y_v} \cap \{t \ge 0 \} \cap \{ \omega=u \} \cr
(x,t) & \mapsto & (\frac{u}{\epsilon}x,\frac{u}{\epsilon}t ) \cr
\end{array}$$
is a homemorphism, since $\widehat{CX_u} \cap \widehat{Y_v} \cap \{t \ge 0 \}$ is conic and $\omega ( \lambda(x,t))=
\lambda \omega(x,t)$ for any $\lambda > 0$. Therefore 
$$\chi \left( {\rm Lk} ( \widehat{CX_u} \cap \widehat{Y_v} \cap \{t \ge 0 \} )\right) =
\chi \left( \widehat{CX_u} \cap \widehat{Y_v} \cap \{t \ge 0 \} \cap \{ \omega=u \} \right).$$
If ${\rm dim} Y \le n-1$ then we can put $\Delta_Y= \mathbb{S}^{n-1} \cap (-Y)$. 

If ${\rm dim} Y=n$ then let $Y'$ be the union of the strata of $Y$ of dimension less than or equal to $n-1$. If $v \in (-Y) \setminus (-Y')$, then by Corollary \ref{Corollary4.17}, we have 
$$ \chi \left( {\rm Lk} (\widehat{CX_u} \cap \widehat{Y_v} \cap \{ t \ge 0 \} )\right)= 
1+ \chi_c \left( \widehat{CX_u} \cap \widehat{Y_v} \cap \{ t \ge 0 \} \cap \{ \omega=\epsilon, t < \frac{2}{a}\epsilon \}
\right),$$
$$ \chi \left( {\rm Lk} (\widehat{X} \cap \widehat{Y_v} \cap \{ t \ge 0 \} )\right)= 
1+ \chi_c \left( \widehat{X} \cap \widehat{Y_v} \cap \{ t \ge 0 \} \cap \{ \omega =\epsilon, t < \frac{2}{a}\epsilon \}
\right),$$
for $0< \epsilon \ll 1$ and where $a$ is such that $\widehat{Y'_v} \cap \{ t \ge 0 \} \subset 
\{ (x,t) \ \vert \ \omega (x) \ge a t \}$ in a neighborhood of $(0,0)$. 

As above let us choose $\epsilon_1 \ge 0$ such that for $0< u \le \epsilon_1$, 
$$\chi \left( {\rm Lk} (\widehat{X} \cap \widehat{Y_v} \cap \{ t \ge 0 \} )\right)= 1+
\chi_c \left( \widehat{X} \cap \widehat{Y_v} \cap \{ t \ge 0 \} \cap \{ \omega =u, t < \frac{2}{a} u \}
\right),$$
and
$$\displaylines{
\qquad \chi_c \left( \widehat{X} \cap \widehat{Y_v} \cap \{ t \ge 0 \} \cap \{ \omega =u, t < \frac{2}{a} u \}
\right)= \hfill \cr
\hfill \chi_c \left( \widehat{CX_u} \cap \widehat{Y_v} \cap \{ t \ge 0 \} \cap \{ \omega =u, t < \frac{2}{a} u \}
\right). \qquad \cr
}$$
Using the same homeomorphism as above, we can conclude that 
$$\chi \left( {\rm Lk} ( \widehat{CX_u} \cap \widehat{Y_v} \cap \{t \ge 0 \} ) \right) = 1 +
\chi_c \left( \widehat{CX_u} \cap \widehat{Y_v} \cap \{ t \ge 0 \} \cap \{ \omega =u, t < \frac{2}{a} u \} \right).$$
So if ${\rm dim} Y =n$, we put $\Delta_Y = \mathbb{S}^{n-1} \cap (-Y')$. \endproof

\begin{proposition}\label{Proposition8.2}
If $Y \subset \mathbb{R}^n$ is a closed conic definable set, then for any germ of closed definable set $(X,0) \subset 
(\mathbb{R}^n,0)$, the following principal kinematic formula holds:
$$ \sigma(X,Y,0)= \sum_{i=0}^n \Lambda_i^{\rm lim} (X,Y,0) \cdot \sigma_{n-i}(Y,0).$$
\end{proposition}
\proof By Lemma \ref{Lemma4.10}, there exists a definable subset $\Sigma_{X,Y} \subset SO(n)$ of positive codimension such that for $\gamma \notin \Sigma_{X,Y}$, $X$ and $\gamma Y$ satisfy Condition (1). Let us fix $\gamma \notin \Sigma_{X,Y}$. By Lemma \ref{Lemma4.11}, there exists a definable subset $\Gamma_{X,\gamma Y} \subset \mathbb{S}^{n-1}$ of positive codimension such that for $v \notin \Gamma_{X,\gamma Y}$, $v$ satisfies Condition (2). Let us choose $v \notin \Gamma_{X,\gamma Y}$. By Lemma \ref{Lemma4.2}, the function $t : \widehat{X} \cap \widehat{(\gamma Y)_v} \to \mathbb{R}$ has an isolated stratified critical point at $(0,0)$. 

Applying Lemma 3.1 in \cite{DutertreJofSingNonIsol} and Lemma \ref{Lemma4.4}, we obtain that
$$\lim_{\epsilon \to 0} \lim_{\delta \to 0^+} 
\chi \left( X \cap (\gamma Y + \delta v) \cap \mathbb{B}_\epsilon^n \right)= \chi \left( {\rm Lk} ( \widehat{X} \cap \widehat{(\gamma Y)_v} \cap \{ t \ge 0 \} ) \right),$$
and 
$$\sigma(X,Y,0)=\frac{1}{s_{n-1}^2} \int_{SO(n) \times \mathbb{S}^{n-1}} \chi \left( {\rm Lk} ( \widehat{X} \cap \widehat{(\gamma Y)_v} \cap \{ t \ge 0 \} )\right) d\gamma dv.$$ 
Of course the same equality is true if we replace $X$ with $CX_u$. By Proposition \ref{Proposition7.2} for $0< u \le \epsilon_0$, we have
$$ \sigma(CX_u,Y,0)= \sum_{i=0}^n  \frac{\Lambda_i(CX_u,CX_u \cap \mathbb{B}^n)}{b_i } \cdot \sigma_{n-i}(Y,0).$$
By Lemma \ref{Lemma8.1}, for $\gamma \in SO(n)$ and $v \notin \Delta_{\gamma Y}$,
$$ \lim_{u \to 0} \chi \left( {\rm Lk} ( \widehat{CX_u} \cap \widehat{Y_v} \cap \{t \ge 0 \} )\right) =
\chi \left( {\rm Lk} ( \widehat{X} \cap \widehat{Y_v} \cap \{t \ge 0 \}) \right).$$ 
Hence, by Hardt's theorem and Lebesgue's theorem,
$$\lim_{u \to 0} \sigma(CX_u,Y,0) =\sigma(X,Y,0).$$
Moreover using Proposition 3.6 in \cite{DutertreGeoDedicata2012}, for $i \in \{1,\ldots,n-2\}$, we have
$$\frac{\Lambda_i(CX_u,CX_u \cap \mathbb{B}^n)}{b_i } =-\frac{1}{2g_n^{n-i-1}}  \int_{G_n^{n-i-1}} \chi (CX_u \cap \mathbb{S}^{n-1} \cap H) dH $$ $$+ \frac{1}{2g_n^{n-i+1}}  \int_{G_n^{n-i+1}} \chi (CX_u \cap \mathbb{S}^{n-1} \cap H) dH.$$
But 
$$\chi (CX_u \cap \mathbb{S}^{n-1} \cap H )= \chi (CX_u \cap \mathbb{S}_u^{n-1} \cap H )= \chi (X \cap \mathbb{S}_u^{n-1} \cap H ),$$
for $H \in G_n^{n-i+1}$ or $H \in G_n^{n-i-1}$, and $\lim_{u \to 0} \chi (CX_u \cap \mathbb{S}^{n-1} \cap H)=\chi \left( {\rm Lk} (X \cap H) \right)$. Passing to the limit as $u \to 0$ and applying Theorem 5.1 in \cite{DutertreGeoDedicata2012}, we get that 
$$\lim_{u \to 0} \frac{\Lambda_i(CX_u,CX_u \cap \mathbb{B}^n)}{b_i} =
\lim_{\epsilon \to 0} \frac{\Lambda_i(X,X\cap \mathbb{B}_\epsilon^n)}{b_i \epsilon^i} .$$
The same proof works for $i=n-1$ and $i=n$. Combining all these equalities, we get the result. \endproof

{\bf Second step:} $X \subset \mathbb{R}^n$ is a closed conic definable set, $(Y,0) \subset (\mathbb{R}^n,0)$ is a germ of closed definable set.

\begin{corollary}\label{Corollary8.3}
Let $X \subset \mathbb{R}^n$ be a closed conic definable set. For any germ of closed definable set $(Y,0) \subset (\mathbb{R}^n,0)$, the following principal kinematic formula holds:
$$\sigma(X,Y,0)= \sum_{i=0}^n \Lambda_i^{\rm lim}(X,0)\cdot \sigma_{n-i}(Y,0).$$
\end{corollary}
\proof By Proposition \ref{Proposition4.20}, we know that $\sigma(X,Y,0)=\sigma(Y,X,0)$, and by Lemma \ref{Lemma4.23} that $\sum_{i=0}^n \Lambda_i^{\rm lim}(X,0)\cdot \sigma_{n-i}(Y,0)=\sum_{i=0}^n \Lambda_i^{\rm lim}(Y,0)\cdot \sigma_{n-i}(X,0)$. Then we apply the previous proposition.
\endproof

{\bf Third step:} $(X,0) \subset (\mathbb{R}^n,0)$ and $(Y,0) \subset (\mathbb{R}^n,0)$ are germs of closed definable sets. 

We assume that $X$ and $Y$ are included in an open neighborhood $U$ of $0$. Let $\epsilon_0 >0$ be such that $\mathbb{B}_{\epsilon_0}^n \subset U$. We set 
$$\Gamma_X = \left\{ (\frac{1}{u}x,u) \ \vert \ x \in X, u \in ]0,\epsilon_0[ \right\} \subset \mathbb{R}^{n+1}$$
and
$$\Gamma_Y = \left\{ (\frac{1}{u}y,u) \ \vert \ y \in Y, u \in ]0,\epsilon_0[ \right\} \subset \mathbb{R}^{n+1}.$$
We recall that the tangent cones of $X$ and $Y$ are $C_0 X = \overline{\Gamma_X} \cap \mathbb{R}^n \times \{0\}$ and $C_0 Y = \overline{\Gamma_Y} \cap \mathbb{R}^n \times \{0\}$. We will now define two tangent cones associated with $\widehat{X}$ and 
$\widehat{Y_v}$ and will relate them to $C_0 X$ and $C_0 Y$. Let 
$$\widehat{\Gamma_X} = \left\{ (\frac{1}{u}x,t,u) \ \vert \ x \in X, u \in ]0,\epsilon_0[ \right\} \subset \mathbb{R}^{n+2}.$$
The following lemma is easy to prove.
\begin{lemma}\label{Lemma8.4}
A point $(x,t)$ belongs to $\overline{\widehat{\Gamma_X}} \cap (\mathbb{R}^{n+1} \times \{0\})$ if and only if there is a sequence of points $(x_n,t_n)_{n \in \mathbb{N}}$ in $\widehat{X}$ and a sequence of positive real numbers $(u_n)_{n \in \mathbb{N}}$ such that $u_n \to 0$ and $(\frac{x_n}{u_n},t_n) \to (x,t)$.
\end{lemma}

\begin{corollary}\label{Corollary8.5}
We have $\widehat{C_0 X}= \overline{\widehat{\Gamma_X}} \cap (\mathbb{R}^{n+1} \times \{0\})$.
\end{corollary}
\proof If $(x,t) \in \widehat{C_0 X}$ then there is a sequence of points $(x_n)_{n \in \mathbb{N}}$ in $X$ and a sequence of positive real numbers $(u_n)_{n \in \mathbb{N}}$ such that $u_n \to 0$ and $\frac{x_n}{u_n} \to x$. Applying the previous lemma to the sequences $(x_n,t)$ and $(u_n)$, we see that $(x,t) \in \overline{\widehat{\Gamma_X}} \cap (\mathbb{R}^{n+1} \times \{0\})$. 

Conversely if $(x,t) \in \overline{\widehat{\Gamma_X}} \cap (\mathbb{R}^{n+1} \times \{0\})$, then there is a sequence of points $(x_n,t_n)_{n \in \mathbb{N}}$ in $X \times \mathbb{R}$ and a sequence of positive real numbers $(u_n)_{n \in \mathbb{N}}$ such that $u_n \to 0$ and $(\frac{x_n}{u_n},t_n) \to (x,t)$. This implies that $x \in C_0 X$ and so that $(x,t) \in \widehat{C_0 X}$. \endproof

Let $v \in \mathbb{S}^{n-1}$ and let 
$$ \widehat{(\Gamma_Y)_v}  =
\left\{ (\frac{1}{u}y,\frac{1}{u}t,u ) \ \vert \ (y,t) \in \widehat{Y_v}, u \in ]0,\epsilon_0[ \right\} \subset \mathbb{R}^{n+2}.$$
\begin{lemma}\label{Lemma8.6}
A point $(y,t)$ belongs to $\overline{\widehat{(\Gamma_Y)_v}} \cap (\mathbb{R}^{n+1} \times \{0\})$ if and only if there is a sequence of points $(y_n,t_n)_{n\in \mathbb{N}}$ in $\widehat{Y_v}$ and a sequence of positive real numbers $(u_n)_{n \in \mathbb{N}}$ such that $u_n \to 0$ and $(\frac{y_n}{u_n},\frac{t_n}{u_n}) \to (y,t)$.
\end{lemma}

\begin{corollary}\label{Corollary8.7}
We have $\widehat{(C_0 Y)_v} = \overline{\widehat{(\Gamma_Y)_v}} \cap (\mathbb{R}^{n+1} \times \{0\})$.
\end{corollary}
\proof If $(y,t) \in \widehat{(C_0 Y)_v}$ then there is a sequence of points $(y_n)_{n\in \mathbb{N}}$ in $Y$ and a sequence of positive real numbers $(u_n)_{n \in \mathbb{N}}$ such that $u_n \to 0$ and $\frac{y_n}{u_n} \to y-tv$. For $n \in \mathbb{N}$, $(y_n+u_ntv, u_ntv)$ is in $\widehat{Y_v}$ and $(\frac{y_n+u_ntv}{u_n}, \frac{u_nt}{u_n})$ tends to $(y,t)$. Therefore $(y,t)$ is in $\overline{\widehat{(\Gamma_Y)_v}} \cap (\mathbb{R}^{n+1} \times \{0\})$.

Conversely if $(y,t)$ is in $\overline{\widehat{(\Gamma_Y)_v}} \cap (\mathbb{R}^{n+1} \times \{0\})$, then there is a sequence of points $(y_n,t_n)_{n \in \mathbb{N}}$ in $\widehat{Y_v}$ and a sequence of positive real numbers $(u_n)_{n \in \mathbb{N}}$ such that $u_n \to 0$ and $(\frac{y_n}{u_n},\frac{t_n}{u_n}) \to (y,t)$. Then $y_n-t_nv \in Y$ and $\frac{y_n-t_nv}{u_n}$ tends to $y-tv$. So $y-tv$ belongs to $C_0 Y$. \endproof

We note that $C_0 X = \widehat{C_0 X} \cap (\mathbb{R}^n \times \{0\})$, $C_0 Y= \widehat{(C_0 Y)_v} \cap (\mathbb{R}^n \times \{0\})$ and that $\widehat{C_0 X}$ and $\widehat{(C_0 Y)_v}$ are closed conic definable sets. 

Let us assume that $X$ is equipped with a Whitney stratification $\mathcal{S}=\{S_i\}_{i=0}^l$ with $S_0=\{0\}$ and $0 \in \overline{S_i}$ for $i=1,\ldots,l$. We set 
$$\Gamma_{S_i}= \left\{ (\frac{x}{u},u) \ \vert \ 
x\in S_i, u \in ]0,\epsilon_0[ \right\} \subset \mathbb{R}^{n+1}$$
for $i=0,\ldots,s$.
\begin{lemma}\label{Lemma8.8}
The partition $\Gamma_X = \cup_{i=0}^s \Gamma_{S_i}$ is a Whitney stratification of $\Gamma_X$.
\end{lemma}
\proof The partition $\cup_{i=0}^s S_i \times ]0,\epsilon_0[$ gives a Whitney stratification of $X \times ]0,\epsilon_0[$. Let $\phi$ be the diffeomorphism $\phi(x,u)= (\frac{1}{u}x,u)$. We have $\Gamma_{S_i}=\phi (S_i \times ]0,\epsilon_0[)$ for $i=0,\ldots,s$ and $\Gamma_X = \phi (X \times ]0,\epsilon_0[)$. This gives the result for Whitney's conditions are $C^1$-invariant. \endproof

We can equip $C_0 X$ with a definable stratification $(\Sigma_k)_{k=0}^{l'}$ where $\Sigma_0=\{0\}$ and $\Sigma_k$ is conic. This is possible for example by considering a Whitney stratification of $C_0 X \cap \mathbb{S}^{n-1}$ and extending it to $C_0 X$ using the conic structure.
\begin{lemma}\label{Lemma8.9}
Let $\Gamma_S$ be a stratum of $\Gamma_X$ and let $\Sigma$ be a stratum of $C_0 X$ such that $\Sigma \subset \overline{\Gamma_S} \setminus \Gamma_S$. The set of points $x$ in $\Sigma$ such that the Thom $(a_u)$-condition is not satisfied at $x$ for the pair $(\Gamma_S,\Sigma)$ is a conic definable set of positive codimension.
\end{lemma}
\proof By \cite{Bekka} and \cite{TaLeLoi}, we already know that this set is definable of positive codimension in $\Sigma$. If $x$ is in this set, then there exists a sequence of points $p_n=(y_n,u_n)$ in $\Gamma_S$ such that $p_n \to (x,0)$ and $T_x \Sigma \nsubseteq \lim_{n \to +\infty} T_{p_n} (\Gamma_S \cap \{u=u_n\} )$. Let $\lambda >0$, by the conic structure of $\Sigma$, $T_{\lambda x} \Sigma =T_x \Sigma$. Moreover $q_n=(\lambda y_n, \frac{u_n}{\lambda}) \in \Gamma_S$ and $T_{p_n} (\Gamma_S \cap \{u=u_n\})= T_{q_n} (\Gamma_S \cap \{u=\frac{u_n}{\lambda}\})$. \endproof

Since the Thom $(a_u)$-condition is stratifying and taking a refinement if necessary, we can assume that the Thom $(a_u)$-condition is satisfied for any pair of strata $(\Gamma_S,\Sigma)$ (see \cite{LojStaWac} for the argument).
This induces a Whitney stratification of $\overline{\widehat{\Gamma_X}}$ compatible with $\widehat{C_0 X} \times \{0\}$. Namely if $S' \subset \mathbb{R}^{n+1}$ is a stratum $\overline{\Gamma_X}$ then 
$\widehat{S'}=\{ (x,t,u) \ \vert \ (x,u) \in S' \}$ is a stratum of $\overline{\widehat{\Gamma_X}}$.
\begin{lemma}\label{Lemma8.10}
This induced stratification of $\overline{\widehat{\Gamma_X}}$ satisfies the Thom $(a_u)$-condition.
\end{lemma}
\proof Let $(x_n,t_n,u_n)_{n\in \mathbb{N}}$ be a sequence of points in $\overline{\widehat{\Gamma_X}}$ that tends to $(x,t,0)$. We can assume that $(x_n,t_n,u_n)$ lies in a stratum $\widehat{S'_1}=\{ (x',t',u') \ \vert \ (x',u') \in S'_1\}$ and that $(x,t,0) \in \widehat{S'_0}$, where $\widehat{S'_0}=\{ (x',t',0) \ \vert \ (x',0) \in S'_0\}$. Since the pair $(S'_1,S'_0)$ satisfies the Thom $(a_u)$-condition, $T_{(x,0)} S'_0 \subset \lim_{n \to +\infty} T_{(x_n,u_n)} (S'_1 \cap \{u=u_n\})$. But 
$$T_{(x_n,t_n,u_n)} \widehat{S'_1}=\left\{ (\nu,\tau,\xi) \ \vert \ (\nu,\xi) \in T_{(x_n,u_n)} S'_1 \right\},$$
and 
$$T_{(x,t,0)} \widehat{S'_0}=\left\{ (\nu,\tau,0) \ \vert \ (\nu,0) \in T_{(x,0)} S'_0 \right\}.$$
It is straightforward to conclude using the fact that $T_{(x_n,u_n)} (S'_1 \cap \{u=u_n\})=
(T_{(x_n,u_n)} S'_1 )\cap \{u=0\}$ and $T_{(x_n,t_n,u_n)} (\widehat{S'_1} \cap \{u=u_n\}) =(T_{(x_n,t_n,u_n)} \widehat{S'_1} )\cap \{u=0 \}$ if $n$ is sufficiently big. \endproof
Similarly we can equip $\overline{\Gamma_Y}$ with a Whitney definable stratification compatible with $C_0 Y \times \{0\}$, that satisfies the Thom $(a_u)$-condition and such that the strata of $C_0 Y$ are conic. This induces a Whitney stratification of $\overline{\widehat{(\Gamma_Y)_v}}$ compatible with $\widehat{(C_0 Y)_v} \times \{0\}$. Namely if $T' \subset \mathbb{R}^{n+1}$ is a stratum of $\overline{\Gamma_Y}$ then  $\widehat{T'_v}= \{ (y,t,u) \ \vert \ (y-tv,u)  \in T' \}$ is a stratum of $\overline{\widehat{(\Gamma_Y)_v}}$ (see Section 4).
\begin{lemma}\label{Lemma8.11}
This induced stratification of $\overline{\widehat{(\Gamma_Y)_v}}$ satisfies the Thom $(a_u)$-condition.
\end{lemma}
\proof The proof is the same as in the previous lemma, taking into account the following remark: if $\widehat{T'_v}=\{(y,t,u) \ \vert \ (y-tv,u) \in T' \}$ is a stratum of $\overline{\widehat{(\Gamma_y)_v}}$ then 
$$T_{(y,t,u)} \widehat{T'_v} =\{ (\nu,\tau,\xi) \ \vert \ (\nu-\tau v ,\xi) \in T_{(y-tv,u)} T' \}.$$ \endproof

For $0<u \le \epsilon_0$, we set $X_u = X \cap \mathbb{S}_u^{n-1}$ and we denote by $CX_u$ the cone over $X_u$, i.e.
$$CX_u =\{ x \in \mathbb{R}^n \ \vert \ \exists \lambda \in \mathbb{R}^+ \hbox{ and }
z \in X_u \hbox{ such that } x=\lambda z \}.$$

\begin{lemma}\label{Lemma8.12} Let us assume that $C_0 X$ and $C_0 Y$ satisfy Condition (1).
There exists a definable subset $\Delta_{X,Y} \subset \mathbb{S}^{n-1}$ of positive codimension such that for $v \notin \Delta_{X,Y}$,
$$ \lim_{u \to 0} \chi \left( {\rm Lk} (\widehat{CX_u} \cap \widehat{Y_v} \cap \{ t\ge 0\}) \right)=
\chi \left( {\rm Lk} (\widehat{X} \cap \widehat{Y_v} \cap \{ t\ge 0\}) \right).$$
\end{lemma}
\proof  As in Lemma \ref{Lemma8.1}, we have to prove that 
$$ \chi \left( {\rm Lk} (\widehat{CX_u} \cap \widehat{Y_v} \cap \{ t\ge 0\}) \right)= 
\chi \left(\widehat{CX_u} \cap \widehat{Y_v} \cap \{ t\ge 0\}\cap \{ \omega=u \} \right),$$
if $u$ is small enough, taking into account that $\widehat{Y_v}$ is not conic. 

First let us fix $v \notin (-C_0 Y) \cap \mathbb{S}^{n-1}$. There exist $\epsilon_v >0$ and $a>0$ such that 
$$\widehat{Y_v} \cap \{ t\ge 0 \} \cap \mathbb{B}_{\epsilon_v}^{n+1} 
\subset \left\{ (x,t) \ \vert \ \omega (x) \ge a t\right\} \cap \mathbb{B}_{\epsilon_v}^{n+1},$$
which implies that there exists $u_0 >0$ such that 
$$\widehat{Y_v} \cap \{ t\ge 0 \} \cap \{ \omega \le u_0 \} \cap \mathbb{B}_{\epsilon_v}^{n+1} \subset \mathring{\mathbb{B}_{\epsilon_v}^{n+1}}.$$

We have assumed that $C_0 X$ and $C_0 Y$ satisfy Condition (1). This means that two strata $W$ and $W'$ of $C_0 X$ and $C_0 Y$ (different from $\{0\}$) intersect transversally. Since these strata are conic, $W\cap \{\omega=1\}$ and $W'$ intersect transversally as well and so $C_0 X \cap \{\omega=1\}$ and $C_0 Y$ intersect transversally (in the stratified sense). As in Lemma \ref{Lemma4.11}, there exists a  definable subset $\Gamma_{C_0 X \cap \{\omega=1\},C_0 Y} \subset \mathbb{S}^{n-1}$ of positive codimension such that for $v \notin \Gamma_{C_0 X \cap \{\omega=1\},C_0 Y}$, 
$\widehat{C_0 X \cap \{\omega= 1\}}$  and $\widehat{(C_0 Y)_v}$ intersect transversally (in the stratified sense).

We need a first auxiliary lemma.
\begin{lemma}\label{lemma8.13}
If $v \notin \Gamma_{C_0 X \cap \{\omega=1\},C_0 Y}$, then there exists $0<u_1 \le u_0$ such that for $0<u \le u_1$ and for $(x,t) \in \widehat{CX_u} \cap \widehat{Y_v} \cap \{t >0\} \cap \{0 < \omega \le u\} \cap \mathbb{B}_{\epsilon_v}^{n+1}$, the sets $\widehat{CX_u}\cap \{t > 0\} \cap \{\omega=\omega(x)\}$ and $\widehat{Y_v} \cap \{ t > 0\}$ intersect transversally (in the stratified sense) at $(x,t)$.
\end{lemma}
\proof Let us specify the stratifications we are working with. The set $Y$ is equipped with a Whitney stratification $\{T_j\}_{j=0}^m$, which induces a stratification $\{ \widehat{(T_j)_v}\}_{j=0}^m$ of $\widehat{Y_v}$. Hence $\widehat{Y_v}\cap \{t>0\}$ is stratified by $\{\widehat{(T_j)_v} \cap \{t>0\}\}_{j=0}^m$. The set $X$ is equipped with a Whitney stratification $\{S_i\}_{i=0}^l$. Hence for $u$ small, $CX_u$ is stratified by $\{0\} \cup \{ C(S_i \cap \mathbb{S}_u^{n-1}) \}_{i=1}^l$. As above this induces a stratification of $\widehat{CX_u} \cap \{t >0\}$. We note that by the conic structure, the intersection $\widehat{CX_u}\cap \{ t> 0\} \cap \{\omega=\omega(x)\}$ is always transverse (in the stratified sense) and the stratification of $\widehat{CX_u} \cap \{t>0\} \cap \{ \omega=\omega(x)\}$ is clear. 

Assume that the above result is not true. Then we can find a sequence of positive real numbers $(u_n)_{n \in \mathbb{N}}$ that tends to $0$ and a sequence of points $(x_n,t_n)_{n \in \mathbb{N}}$ in $\widehat{CX_{u_n}}\cap \widehat{Y_v} \cap \{t>0\} \cap \{0< \omega \le u_n\} \cap \mathbb{B}_\epsilon^{n+1}$ such that $\widehat{CX_{u_n}} \cap \{\omega  =\omega (x_n) \}$ and $\widehat{Y_v}$ do not intersect transversally at $(x_n,t_n)$.

We can assume that the sequence $(x_n,t_n)_{n \in \mathbb{N}}$ is included in a unique stratum $\widehat{T_v}$, where $T$ is a stratum of $Y$. Moreover we can assume that there is a stratum $S\not= \{0\}$ of $X$ such that for each $n\in \mathbb{N}$, $(x_n,t_n) \in \widehat{CS_{u_n}}$. Since $\frac{t_n}{\omega(x_n)} \le \frac{1}{a}$, taking a subsequence if necessary, we can assume that $(\frac{x_n}{\omega(x_n)},\frac{t_n}{\omega(x_n)})$ tends to $(x,t)$. Let $w_n=u_n \frac{x_n}{\omega(x_n)}$, then $w_n \in CS_{u_n} \cap \{\omega=u_n\}$ and so $w_n \in S \subset X$. Therefore $(w_n, \frac{t_n}{\omega(x_n)})$ belongs to $\widehat{X}$ and $(\frac{1}{u_n}w_n,\frac{t_n}{\omega(x_n)})$ tends to $(x,t)$. By Lemma 8.4 and Corollary 8.5, this implies that $(x,t)$ is in $\widehat{C_0 X}$. Moreover since $\omega (\frac{w_n}{u_n})=1$, $(x,t)$ belongs to $\widehat{C_0 X} \cap \{ \omega=1\}$. On the other hand, $(x_n,t_n) \in \widehat{Y_v}$ and so by Lemma 8.6 and Corollary 8.7, $(x,t)$ belongs to $\widehat{(C_0Y)_v}$. 

The points $p_n:=(\frac{1}{u_n}w_n,\frac{t_n}{\omega (x_n)},u_n)$ belong to the stratum $\widehat{\Gamma_S}$ of $\widehat{\Gamma_X}$ and the point $(x,t)$ belongs to a stratum $\widehat{\Sigma}$ of $\widehat{C_0 X}$. By the Thom $(a_u)$-condition, we have $T_{(x,t)} \widehat{\Sigma} \subset \lim_{n \to +\infty} T_{p_n} (\widehat{\Gamma_S} \cap \{u=u_n\})$.

Since $\widehat{\Sigma}$ is $\Sigma \times \mathbb{R}$ and $\Sigma$ is conic, $\widehat{\Sigma}$ intersects $\{\omega=1\}$ transversally. By the Thom $(a_u)$-condition, $\{ \omega = 1\}$ intersects $\widehat{\Gamma_S}\cap \{u=u_n\}$ transversally for $n$ big enough, and so 
$$T_{(x,t)}( \widehat{\Sigma}\cap \{\omega=1\}) \subset \lim_{n \to +\infty} T_{p_n}( \widehat{\Gamma_S}\cap \{\omega=1\}  \cap \{u=u_n\}).$$
But $T_{p_n} (\widehat{\Gamma_S}\cap \{\omega=1\}  \cap \{u=u_n\} )= T_{(w_n,\frac{t_n}{\omega (x_n)})} (\widehat{S} \cap \{\omega=u_n\})$ and so 
$$ T_{(x,t)} (\widehat{\Sigma}\cap \{\omega=1\}) \subset \lim_{n \to +\infty} T_{(w_n,\frac{t_n}{\omega (x_n)})} (\widehat{S} \cap \{\omega=u_n\}).$$
We note that
$$ T_{(w_n,\frac{t_n}{\omega(x_n)})} (\widehat{S} \cap \{\omega=u_n\} )=
T_{(x_n,t_n)} (\widehat{CS_{u_n}} \cap \{ \omega = \omega (x_n) \}),$$
by the conic structure of $CS_{u_n}$. 

The points $q_n:= (\frac{x_n}{\omega (x_n)}, \frac{t_n}{\omega(x_n)},\omega (x_n) )$ are in the stratum $\widehat{(\Gamma_T)_v}$ of $\widehat{(\Gamma_Y)_v}$ and the point $(x,t)$ is in a stratum $\widehat{\Sigma'_v}$ of $\widehat{(C_0Y)_v}$. By the Thom $(a_u)$-condition, we have $T_{(x,t)} \widehat{\Sigma'_v} \subset \lim_{n \to + \infty} T_{q_n} ( \widehat{(\Gamma_T)_v} \cap \{u = \omega(x_n) \})$. But $T_{q_n} (\widehat{(\Gamma_T)_v} \cap \{ u = \omega(x_n) \}) = T_{(x_n,t_n)} \widehat{T_v}$ and so 
$T_{(x,t)} \widehat{\Sigma'_v} \subset \lim_{n \to +\infty} T_{(x_n,t_n)} \widehat{T_v}$. 

Since $v \notin \Gamma_{C_0 X \cap \{\omega=1\},C_0 Y}$, 
$\widehat{C_0 X \cap \{\omega= 1\}}$  and $\widehat{(C_0 Y)_v}$ intersect transversally (in the stratified sense). But $\widehat{C_0 X \cap \{\omega= 1\}} = \widehat{C_0 X} \cap \{ \omega=1\}$, and we conclude that 
$$T_{(x,t)} (\widehat{\Sigma}\cap \{\omega=1\}) + T_{(x,t)} \widehat{\Sigma'_v}=\mathbb{R}^{n+1}.$$
Therefore 
$$\lim_{n \to +\infty} T_{(x_n,t_n)} (\widehat{CS_{u_n}} \cap \{ \omega=\omega(x_n) \} )+
\lim_{n \to +\infty} T_{(x_n,t_n)} \widehat{T_v} =\mathbb{R}^{n+1},$$
and so, for $n$ big enough 
$$T_{(x_n,t_n)} (\widehat{CS_{u_n}} \cap \{ \omega=\omega(x_n) \}) + T_{(x_n,t_n)}\widehat{T_v}=\mathbb{R}^{n+1}.$$
This contradicts the construction of the sequence $(x_n,t_n)$ and ends the proof of this auxiliary lemma. \endproof

Similarly the following second auxiliary lemma holds.
\begin{lemma}\label{lemma8.14}
There exists $0<u_2 \le u_0$ such that for $0<u \le u_2$ and for $x \in CX_u \cap Y  \cap \{0 < \omega \le u\} \cap \mathbb{B}_{\epsilon_v}^{n}$, the sets $CX_u \cap \{\omega=\omega(x)\}$ and $Y$ intersect transversally (in the stratified sense) in $\mathbb{R}^n$ at $x$.
\end{lemma}

Let us choose $u>0$ such that $u \le {\rm min}\{\epsilon_1,u_1,u_2\}$, where $\epsilon_1$ is such that for $0< u \le \epsilon_1$,
$$\chi \left( {\rm Lk} ( \widehat{X} \cap \widehat{Y_v} \cap \{t \ge 0 \} ) \right) =
\chi \left(\widehat{X} \cap \widehat{Y_v} \cap \{t \ge 0 \} \cap \{ \omega=u \} \right).$$  Then for $(x,t) \in \widehat{CX_u}\cap \widehat{Y_v}\cap \{ 0< \omega \le u\} \cap \{ t> 0 \} \cap \mathbb{B}_{\epsilon_v}^{n+1}$, $\widehat{CX_u}\cap \{ \omega=\omega(x)\} \cap \{ t>0\}$ and $\widehat{Y_v} \cap \{ t>0\}$ intersect transversally (in the stratified sense) at $(x,t)$. This implies that $\widehat{CX_u} \cap \{ t>0\}$ and $\widehat{Y_v} \cap \{ t>0\}$ intersect transversally at $(x,t)$ and that $\{\omega=\omega(x)\}$ intersects $\widehat{CX_u}\cap \widehat{Y_v} \cap \{ t>0\}$ transversally at $(x,t)$, and so $(x,t)$ is not a stratified critical point of $\omega_{\vert \widehat{CX_u} \cap \widehat{Y_v} \cap \{t>0\}}$. Similarly if $(x,0) \in \widehat{CX_u}\cap \widehat{Y_v} \cap \{ 0< \omega \le u\} \cap \mathbb{B}_{\epsilon_v}^{n+1}$, then $(x,0)$ is not a stratified critical point of $\omega_{\vert \widehat{CX_u} \cap \widehat{Y_v} \cap \{t=0\}}$. Hence we conclude that $\omega : \widehat{CX_u}\cap \widehat{Y_v}\cap \{ 0< \omega \le u\} \cap \{ t \ge 0 \} \cap \mathbb{B}_{\epsilon_v}^{n+1} \to \mathbb{R}$ is a stratified submersion and so that 
$$\chi \left( {\rm Lk} (\widehat{CX_u} \cap \widehat{Y_v} \cap \{ t \ge 0 \} )\right)=
\chi \left( \widehat{CX_u} \cap \widehat{Y_v} \cap \{ t \ge 0 \} \cap \{\omega=u\} \right).$$
Therefore if ${\rm dim}C_0 Y \le n-1$, we can take $$\Delta_{X,Y} =\left ((-C_0Y) \cap \mathbb{S}^{n-1}\right) \cup \Gamma_{C_0 X\cap \{\omega=1\},C_0 Y}.$$ 

If ${\rm dim} C_0 Y=n$ then ${\rm dim} Y=n$. Let $Y'$ be the union of the strata of $Y$ of dimension less than or equal to $n-1$. If $v \in (-C_0 Y) \setminus (-C_0 Y')$, we know that
$$\chi \left( {\rm Lk} (\widehat{CX_u} \cap \widehat{Y_v} \cap \{t \ge 0\}) \right)= 1 + 
\chi_c \left( \widehat{CX_u} \cap \widehat{Y_v} \cap \{t \ge 0\} \cap \{\omega =\epsilon, t < \frac{2}{a} \epsilon\}
\right),$$
$$\chi \left( {\rm Lk} (\widehat{X} \cap \widehat{Y_v} \cap \{t \ge 0\}) \right)= 1 + 
\chi_c \left( \widehat{X} \cap \widehat{Y_v} \cap \{t \ge 0\} \cap \{\omega =\epsilon, t < \frac{2}{a} \epsilon\}
\right),$$
for $0<\epsilon \ll 1$ and where $a>0$ is such that $\widehat{Y'_v} \cap \{ t\ge 0\} \subset 
\{ (x,t)\ \vert \ \omega (x) \ge at \}$ in a neighborhood of $(0,0)$. 

As above let us choose $\epsilon_1 >0$ such that for $0< u \le \epsilon_1$,
$$\chi \left( {\rm Lk} (\widehat{X} \cap \widehat{Y_v} \cap \{t \ge 0\}) \right)= 1 + 
\chi_c \left( \widehat{X} \cap \widehat{Y_v} \cap \{t \ge 0\} \cap \{\omega =u, t < \frac{2}{a} u\}
\right).$$
Since
$$\displaylines{ \qquad  \chi_c \left( \widehat{X} \cap \widehat{Y_v} \cap \{t \ge 0\} \cap \{\omega =u, t < \frac{2}{a} u \} \right) \hfill \cr
\hfill = \chi_c \left( \widehat{CX_u} \cap \widehat{Y_v} \cap \{t \ge 0\} \cap \{\omega =u, t < \frac{2}{a} u \} \right), \qquad \cr
}$$
we have to prove that 
$$\chi \left( {\rm Lk} (\widehat{CX_u} \cap \widehat{Y_v} \cap \{t \ge 0\}) \right)= 1 + 
\chi_c \left( \widehat{CX_u} \cap \widehat{Y_v} \cap \{t \ge 0\} \cap \{\omega =u, t < \frac{2}{a} u\} \right),$$
if $u$ is small enough. By the previous case, we know that for $u$ small enough and for $v \notin \Gamma_{C_0X \cap \{\omega=1\},C_0 Y'}$, $$\omega : \widehat{CX_u}\cap \widehat{Y'_v}\cap \{ 0< \omega \le u\} \cap \{ t \ge 0 \} \cap \mathbb{B}_\epsilon^{n+1} \to \mathbb{R}$$ is a stratified submersion, for an appropriate $\epsilon >0$. Since the strata of $\widehat{Y_v}\setminus \widehat{Y'_v}$ have dimension $n+1$ and the strata of $Y\setminus Y'$ have dimension $n$, $$\omega :  \widehat{CX_u}\cap \widehat{Y_v}\cap \{ 0< \omega \le u\} \cap \{ t \ge 0 \} \cap \mathbb{B}_\epsilon^{n+1} \to \mathbb{R}$$ is a stratified submersion by the conic structure of $CX_u$. For the same reason and because $\widehat{Y'_v}\cap \{ t\ge 0\} \cap \{ t = \frac{2}{a}\omega(x) \}=\{(0,0)\}$, we see that 
$$\omega :  \widehat{CX_u}\cap \widehat{Y_v}\cap \{ 0< \omega \le u\} \cap \{ t=\frac{2}{a} \omega (x) \}  \to \mathbb{R}$$ is a stratified submersion. Hence for $0< \epsilon \le u$,
$$\displaylines{
\qquad \chi_c \left( \widehat{CX_u} \cap \widehat{Y_v} \cap \{t \ge 0\} \cap \{\omega =u, t < \frac{2}{a} u\} \right) = \hfill \cr
\hfill \chi_c \left( \widehat{CX_u} \cap \widehat{Y_v} \cap \{t \ge 0\} \cap \{\omega =\epsilon, t < \frac{2}{a} \epsilon \} \right). \qquad \cr
}$$
If ${\rm dim} C_0 Y =n$, we take 
$$\Delta_{X,Y} = \left( (-C_0 Y') \cap \mathbb{S}^{n-1} \right) \cup \Gamma_{C_0 X \cap \{\omega=1\}, C_0 Y' }.$$
\endproof

\begin{theorem}\label{Theorem8.15}
Let $(X,0) \subset (\mathbb{R}^n,0)$ and $(Y,0) \subset (\mathbb{R}^n,0)$ be two germs of closed definable sets. The following principal kinematic formula holds:
$$ \sigma(X,Y,0)= \sum_{i=0}^n \Lambda_i^{\rm lim}(X,0) \cdot \sigma_{n-i}(Y,0).$$
\end{theorem}
\proof By Lemma \ref{Lemma4.10}, there exists a definable subset $\Sigma_{X,Y} \subset SO(n)$ of positive codimension such that for $\gamma \notin \Sigma_{X,Y}$, $X$ and $\gamma Y$ satisfy Condition (1). Let us fix $\gamma \notin \Sigma_{X,Y}$. By Lemma \ref{Lemma4.11}, there exists a definable subset $\Gamma_{X,\gamma Y} \subset \mathbb{S}^{n-1}$ of positive codimension such that for $v \notin \Gamma_{X,\gamma Y}$, $v$ satisfies Condition (2). Let us choose $v \notin \Gamma_{X,\gamma Y}$. By Lemma \ref{Lemma4.2}, the function $t : \widehat{X} \cap \widehat{(\gamma Y)_v} \to \mathbb{R}$ has an isolated stratified critical point at $(0,0)$. 

Applying Lemma 3.1 in \cite{DutertreJofSingNonIsol} and Lemma \ref{Lemma4.4}, we obtain that
$$\lim_{\epsilon \to 0} \lim_{\delta \to 0^+} 
\chi \left( X \cap (\gamma Y + \delta v) \cap \mathbb{B}_\epsilon^n \right)= \chi \left( {\rm Lk} ( \widehat{X} \cap \widehat{(\gamma Y)_v} \cap \{ t \ge 0 \})\right),$$
and so
$$\sigma(X,Y,0) = \frac{1}{s_{n-1}^2} \int_{SO(n) \times \mathbb{S}^{n-1}} \chi \left( {\rm Lk} ( \widehat{X} \cap \widehat{(\gamma Y)_v} \cap \{ t \ge 0 \} )\right) dv d\gamma.$$ 
Of course the same equality is true if we replace $X$ with $CX_u$. By Corollary \ref{Corollary8.3}, for $0< u \le \epsilon_0$, we have
$$ \sigma(CX_u,Y,0)= \sum_{i=0}^n  \frac{\Lambda_i(CX_u,CX_u \cap \mathbb{B}^n)}{b_i } \cdot \sigma_{n-i}(Y,0).$$
Since $C_0 (\gamma Y)=\gamma (C_0 Y)$, by Lemma \ref{Lemma4.10} there exists a definable subset $\Sigma_{C_0X,C_0 Y} \subset SO(n)$ of positive codimension such that for $\gamma \notin \Sigma_{C_0 X,C_0 Y}$, $C_0 X$ and $C_0 (\gamma Y)$ satisfy Condition (1). For $\gamma \notin \Sigma_{C_0 X,C_0 Y}$ and $v \notin \Delta_{X,\gamma Y}$,
$$ \lim_{u \to 0} \chi \left( {\rm Lk} ( \widehat{CX_u} \cap \widehat{(\gamma Y)_v} \cap \{t \ge 0 \} )\right) =
\chi \left( {\rm Lk} ( \widehat{X} \cap \widehat{(\gamma Y)_v} \cap \{t \ge 0 \}) \right).$$ 
Hence, by Hardt's theorem and Lebesgue's theorem,
$$\lim_{u \to 0} \sigma(CX_u,Y,0) =\sigma(X,Y,0).$$
We end the proof as in Proposition \ref{Proposition8.2}. \endproof

Let us specify this kinematic formula when $d+e=n$, $d={\rm dim} X$ and $e={\rm dim} Y$.
We denote by $X^d$ (resp. $Y^e$) the union of the top-dimensional strata of $X$ (resp. $Y$). 
\begin{corollary}\label{CorollaryofTh8.15}
Let $(X,0) \subset (\mathbb{R}^n,0)$ and $(Y,0) \subset (\mathbb{R}^n,0)$ be two germs of closed definable sets such that $d+e=n$, where $d={\rm dim} X$ and $e={\rm dim} Y$.The following formula holds: 
$$ 
\frac{1}{s_{n-1}^2} \int_{SO(n) \times \mathbb{S}^{n-1}} \lim_{\epsilon \to 0} \lim_{\delta \to 0^+}  \# \left( X^d \cap (\gamma Y^e + \delta v) \cap \mathbb{B}_\epsilon^n \right) d\gamma dv=
\Theta_d(X) \cdot \Theta_e(Y).$$
\end{corollary}
\proof For $\gamma$ generic in $SO(n)$ and $v$ generic in $\mathbb{S}^{n-1}$, 
$$\lim_{\epsilon \to 0} \lim_{\delta \to 0^+} 
 \chi \left( X \cap (\gamma Y + \delta v) \cap \mathbb{B}_\epsilon^n \right) =
 \lim_{\epsilon \to 0} \lim_{\delta \to 0^+} 
 \# \left( X^d \cap (\gamma Y^e + \delta v) \cap \mathbb{B}_\epsilon^n \right).$$ \endproof

Let us formulate now the second principal kinematic formula.
\begin{theorem}\label{Theorem8.16}
Let $(X,0) \subset (\mathbb{R}^n,0)$ and $(Y,0) \subset (\mathbb{R}^n,0)$ be two germs of closed definable sets. The following principal kinematic formula holds:
$$\Lambda_0^{\rm lim} (X,Y,0)= \sum_{i=0}^n \Lambda_i^{\rm lim}(X,0) \cdot \Lambda_{n-i}^{\rm lim}(Y,0).$$
\end{theorem}
\proof Let us compute the integrals 
$$\frac{1}{s_{n-1}} \int_{SO(n)} \chi \left( {\rm Lk} (X \cap \gamma Y) \right) d\gamma,$$
and 
$$\frac{1}{s_{n-1}^2} \int_{SO(n)}\int_{\mathbb{S}^{n-1}} \chi \left( {\rm Lk} (X \cap \gamma Y\cap \{u^*=0 \}) \right)  dud\gamma.$$
Let us assume first that $X$ and $Y$ are conic closed definable sets. We have already computed the first integral in the proof of Proposition \ref{Proposition5.1} and we have found that
$$\displaylines{
\qquad \frac{1}{s_{n-1}} \int_{SO(n)} \chi \left( {\rm Lk} (X \cap \gamma Y) \right) d\gamma \hfill \cr
\hfill =\sum_{i=0}^{n-1} \frac{\tilde{\Lambda}_i \left( {\rm Lk} (X), {\rm Lk} (X)\right)}{s_i}
\cdot \frac{1}{g_n^{i+1}}\int_{G_n^{i+1}} \chi \left( {\rm Lk} (Y \cap H) \right)dH,\qquad \cr
}$$
which can be rewritten in the following way:
$$\displaylines{
\qquad \frac{1}{s_{n-1}} \int_{SO(n)} \chi \left( {\rm Lk} (X \cap \gamma Y) \right) d\gamma \hfill \cr
\hfill =\sum_{i=1}^{n} \frac{\Lambda_i( X, X \cap \mathbb{B}^n)}{b_i}
\cdot \frac{1}{g_n^{i}}\int_{G_n^{i}} \chi \left( {\rm Lk} (Y \cap H ) \right)dH.
\qquad \cr
}$$
The same computation applied to $X \cap \{u^*=0\}$ yields
$$\displaylines{
\qquad \frac{1}{s_{n-1}^2} \int_{SO(n)}\int_{\mathbb{S}^{n-1}} \chi \left( {\rm Lk} (X \cap \gamma Y\cap \{u^*=0 \}) \right)  dud\gamma \hfill \cr
\qquad \qquad = \sum_{i=0}^{n-1} \Big( \frac{1}{s_{n-1}} \int_{\mathbb{S}^{n-1}}
 \frac{\tilde{\Lambda}_i \left( {\rm Lk} (X \cap \{u^*=0\}), {\rm Lk} (X \cap \{u^*=0\}) \right)}{s_i} du \hfill \cr
\hfill  \times  \frac{1}{g_n^{i+1}}\int_{G_n^{i+1}} \chi \left( {\rm Lk} (Y \cap H) \right)dH \Big) \qquad \cr
\qquad \qquad = \sum_{i=0}^{n-1}  \Big( \frac{1}{g_n^{n-1}} \int_{G_n^{n-1}}
 \frac{\tilde{\Lambda}_i \left( {\rm Lk} (X \cap L), {\rm Lk} (X \cap L) \right)}{s_i} dL \hfill \cr
\hfill \times \frac{1}{g_n^{i+1}}\int_{G_n^{i+1}} \chi \left( {\rm Lk} (Y \cap H) \right)dH. \Big) \qquad \cr
}$$
Using the notations and normalizations of \cite{BernigBroeckerFourier}, Theorem 4.4, we can write
$$\displaylines{
\qquad \frac{1}{g_n^{n-1}} \int_{G_n^{n-1}}
 \frac{\tilde{\Lambda}_i \left( {\rm Lk} (X \cap L), {\rm Lk} (X \cap L)\right)}{s_i} dL \hfill \cr
 \hfill =\frac{1}{s_{n-1}}\int_{SO(n)} \frac{\tilde{\Lambda}_i \left( {\rm Lk} (X \cap \gamma E), {\rm Lk} (X \cap \gamma E) \right)}{s_i} d\gamma, \qquad \cr
}$$
where $E$ is a $(n-2)$-dimensional unit sphere in $\mathbb{S}^{n-1}$. By the spherical kinematic formula, we find that
$$\frac{1}{s_{n-1}}\int_{SO(n)} \frac{\tilde{\Lambda}_i( {\rm Lk} \left(X \cap \gamma E), {\rm Lk} (X \cap \gamma E)\right)}{s_i} d\gamma = \frac{1}{s_{i+1}}\tilde{\Lambda}_{i+1}({\rm Lk}(X),{\rm Lk}(X))$$ $$= \frac{1}{b_{i+2}} \Lambda_{i+2}(X,X \cap \mathbb{B}^n).$$
Hence we get that
$$\displaylines{
\qquad \frac{1}{s_{n-1}^2} \int_{SO(n)}\int_{\mathbb{S}^{n-1}} \chi \left( {\rm Lk} (X \cap \gamma Y\cap \{u^*=0 \}) \right)  dvd\gamma \hfill \cr
\hfill= \sum_{i=2}^n \frac{\tilde{\Lambda}_{i}(X,X \cap \mathbb{B}^n)}{b_i}\cdot \frac{1}{g_n^{i-1}}
\int_{G_n^{i-1}} \chi \left( {\rm Lk}(Y \cap H) \right) dH. \qquad \cr
}$$
Then we apply this result to $CX_u$ and $CY_u$ where $X_u=X \cap \mathbb{S}^{n-1}_u$ and $Y_u=Y \cap \mathbb{S}^{n-1}_u$, and make $u \to 0$ and obtain
$$\frac{1}{s_{n-1}} \int_{SO(n)} \chi \left( {\rm Lk} (X \cap \gamma Y) \right) d\gamma=
\sum_{i=1}^{n} \lim_{\epsilon \to 0} \frac{\Lambda_i( X, X \cap \mathbb{B}_\epsilon^n)}{b_i\epsilon^i}
\cdot \frac{1}{g_n^{i}}\int_{G_n^{i}} \chi \left( {\rm Lk} (Y \cap H ) \right)dH,$$
and
$$\displaylines{
\qquad \frac{1}{s_{n-1}^2} \int_{SO(n)}\int_{\mathbb{S}^{n-1}} \chi \left( {\rm Lk} (X \cap \gamma Y\cap \{u^*=0 \}) \right)  dud\gamma \hfill \cr
\hfill =\sum_{i=2}^{n} \lim_{\epsilon \to 0} \frac{\Lambda_i( X, X \cap \mathbb{B}_\epsilon^n)}{b_i\epsilon^i}
\cdot \frac{1}{g_n^{i-1}}\int_{G_n^{i-1}} \chi \left( {\rm Lk} (Y \cap H ) \right)dH. \qquad \cr
}$$
Therefore by the relation between $\sigma(X,Y,0)$ and $\Lambda_0(X,Y,0)$
and by Theorem \ref{Theorem8.15}, we get 
$$\displaylines{
\qquad \Lambda_0(X,Y,0)=
\Lambda_0^{\rm lim}(X,0) \sigma_n(Y,0) 
+  \hfill \cr
\qquad \qquad \Lambda_1^{\rm lim}(X,0) \left( \sigma_{n-1}(Y,0) -\frac{1}{2g_n^1}\int_{G_n^1} \chi \left( {\rm Lk} (Y \cap H) \right)
dH \right) \hfill \cr
\hfill + \sum_{i=2}^n \Lambda_i^{\rm lim}(X,0) \cdot A_i ,\qquad \qquad \cr
}$$
where
$$A_i= \sigma_{n-i}(Y,0) -\frac{1}{2g_n^i}\int_{G_n^i} \chi \left( {\rm Lk} (Y \cap H) \right) dH -\frac{1}{2g_n^{i-1}}\int_{G_n^{i-1}} \chi \left( {\rm Lk} (Y \cap H) \right)
dH .
$$
By \cite{DutertreJofSingProcTrot}, Theorem 5.6 and its proof, we have that
$\sigma_n(Y,0)= \Lambda_n^{\rm lim}(Y,0)$
and for $i \ge 1$,
$$\sigma_{n-i} (Y,0)= \frac{1}{2g_n^{i+1}}\int_{G_n^{i+1}} \chi \left( {\rm Lk} (Y \cap H) \right) dH + \frac{1}{2g_n^{i}}\int_{G_n^{i}} \chi \left( {\rm Lk} (Y \cap H) \right) dH.$$
Moreover by \cite{DutertreGeoDedicata2012}, Theorem 5.1, we have that
$$\Lambda_{n-1}^{\rm lim}(Y,0)= \frac{1}{2g_n^{2}}\int_{G_n^{2}} \chi \left( {\rm Lk} (Y \cap H) \right) dH$$
and for $i \ge 2$,
$$\Lambda_{n-i}^{\rm lim}(Y,0) = \frac{1}{2g_n^{i+1}}\int_{G_n^{i+1}} \chi \left( {\rm Lk} (Y \cap H) \right) dH
-\frac{1}{2g_n^{i-1}}\int_{G_n^{i-1}} \chi \left( {\rm Lk} (Y \cap H) \right) dH.$$
These equalities enable us to end the proof. \endproof

For $Y=H$, where $H \in G_n^{n-k}$ and $k \in \{1,\ldots,n\}$, the above kinematic formula writes $$ \sigma(X,H,0)=  \Lambda_k^{\rm lim} (X,0).$$
Hence we recover our Theorem \ref{Theorem3.7}, because for $H \in G_n^{n-k}$, 
$$\beta_0(H)= \frac{1}{s_{n-1}} \int_{\mathbb{S}^{n-1}}
\lim_{\epsilon \rightarrow 0} \lim_{\delta \rightarrow 0^+} \Lambda_0(H_{\delta,v} \cap X, H_{\delta,v} \cap X \cap \mathbb{B}^n_\epsilon) dv,$$
by the co-area formula.

\section{More kinematic formulas}
In view of Theorem \ref{Theorem8.16}, a natural question is to express the following sums
$$ \sum_{i+j=p+n} \Lambda_i (X,0) \cdot \Lambda_j (X,0),$$
for $k=1,\ldots,n$ as the right-hand side of a kinematic formula. The answer is quite simple and explained briefly in the next proposition.
\begin{proposition}\label{Proposition9.1}
Let $(X,0) \subset (\mathbb{R}^n,0)$ and $(Y,0) \subset (\mathbb{R}^n,0)$ be two germs of closed definable sets. For $k \in \{1,\ldots,n\}$, the following kinematic formula holds:
$$\int_{SO(n)} \lim_{\epsilon \to 0}
\frac{\Lambda_k(X \cap \gamma Y,X \cap \gamma Y \cap \mathbb{B}^n_\epsilon) }{b_k \epsilon^k} d \gamma= \sum_{i+j=k+n}\Lambda_i^{\rm lim}(X,0) \cdot \Lambda_{j}^{\rm lim}(Y,0).$$
\end{proposition}
\proof When $X$ and $Y$ are conic, it enough to apply the spherical kinematic formulas. 
The general case can be deduced as we have already done in several previous proofs. \endproof

\begin{corollary}\label{Corollary9.2}
Let $(X,0) \subset (\mathbb{R}^n,0)$ and $(Y,0) \subset (\mathbb{R}^n,0)$ be two germs of closed definable sets. The following principal kinematic formula holds:
$$\int_{SO(n)} \lim_{\epsilon \to 0}
\Lambda_0(X \cap \gamma Y,X \cap \gamma Y \cap \mathbb{B}^n_\epsilon)  d \gamma= 
\sum_{i=0}^n \Lambda_i^{\rm lim} (X,0) \cdot \left( \sum_{j=0}^{n-i} \Lambda_j^{\rm lim} (Y,0) \right).$$
\end{corollary}
\proof Apply Proposition \ref{Proposition9.1} and Corollary \ref{Corollary3.6}. \endproof

\end{document}